\documentclass[a4paper,11pt]{article}         
\usepackage[sc]{mathpazo} 

\usepackage{amsmath}  
\usepackage{amsfonts} 
\usepackage{amsthm}
\usepackage{graphicx}
\usepackage{array} 
\usepackage{enumerate} 
\usepackage[latin1]{inputenc} 
\usepackage[T1]{fontenc} 
\usepackage{url}
\usepackage[usenames]{xcolor}
\definecolor{DarkGreen}{rgb}{0,0.5,0}
\usepackage{overpic}
\usepackage{eurosym}
\usepackage{framed}\colorlet{shadecolor}{black!05!white} 
\usepackage{tikz}
\usepackage[misc]{ifsym}
\usepackage{fancyhdr}
\fancyheadoffset{30pt} 
\fancyfootoffset{30pt} 

\newcommand\unio{\mathop\cup\limits}
\newcommand{\Reals}{\mathbb{R}}  
\newcommand{\Naturals}{\mathbb{N}}  
  
\newcommand{\E}{\operatorname{E}}
\newcommand{\ind}{\mathbf{1}}
\newcommand{\Binom}{\operatorname{Binom}}
\newcommand{\Exp}{\operatorname{Exp}}
\newcommand{\Geom}{\operatorname{Geom}}
\newcommand{\Pois}{\operatorname{Pois}}

\newcommand{\Conc}{\textbf}

\def\unio{\mathop\cup\limits} 
\def\Naturals{\mathbb N} 
\def\Reals{\mathbb R} 
\def\Enters{\mathbb Z} 
  \newcommand{\condprob}[2]{\raise2pt
    \hbox{%
      \mathsurround=0pt$#1$}
    \ \raise-1pt\hbox{\scalebox{1.2}[1.5]{/}}\,%
    \raise-2pt
    \hbox{%
      \mathsurround=0pt$#2$}
  }

\def\E{\operatorname{E}} 
\def\Var{\operatorname{Var}} 
\newdimen\fletxauu
\newdimen\fletxadosu
\newdimen\fletxaudos
\newdimen\fletxadosdos
\newdimen\difprimera
\newdimen\difsegona
\newdimen\diftercera
\newbox\Zero
\newbox\Una
\newbox\Dues
\newbox\Tres
\newbox\Quatre
\newbox\Cinc
\newbox\Sis
\newbox\Set
\newbox\Vuit
\newbox\Nou
\def\mapstofill{$\mathsurround=0pt\mapstochar\mathrel{\mkern-4mu}
                \mathord- \mkern-6mu
                \cleaders\hbox{$\mkern-2mu\mathord-\mkern-2mu$}
                \hfill\mkern-6mu\mathord\rightarrow$}
\def\composiciodeaplicacions: nom: #1 conjunt1: #2 aplicacio1: #3
     conjunt2: #4 aplicacio2: #5 conjunt3: #6 element1: #7
     element2: #8 element3: #9 fi{
    \fletxauu=1.2cm
    \fletxadosu=1.2cm
    \fletxaudos=1.2cm
    \fletxadosdos=1.2cm
    \setbox\Una=\hbox{#1}
    \setbox\Dues=\hbox{#2}
    \setbox\Tres=\hbox{#3}
    \setbox\Quatre=\hbox{#4}
    \setbox\Cinc=\hbox{#5}
    \setbox\Sis=\hbox{#6}
    \setbox\Set=\hbox{#7}
    \setbox\Vuit=\hbox{#8}
    \setbox\Nou=\hbox{#9}
    \difprimera=\wd\Tres
    \advance\difprimera by -0.4cm
    \ifdim\difprimera>0pt\advance\fletxauu by \difprimera
                         \advance\fletxadosu by \difprimera
                         \fi
    \difprimera=\wd\Dues
    \advance\difprimera by -\wd\Set
    \divide\difprimera by 2
    \ifdim\difprimera>0pt\advance\fletxadosu by \difprimera
                  \else\advance\fletxauu by -\difprimera
          \fi
    \difsegona=\wd\Cinc
    \advance\difsegona by -0.4cm
    \ifdim\difsegona>0pt\advance\fletxaudos by \difsegona
                        \advance\fletxadosdos by \difsegona
                        \fi
    \difsegona=\wd\Quatre
    \advance\difsegona by -\wd\Vuit
    \divide\difsegona by 2
    \ifdim\difsegona>0pt\advance\fletxadosu by \difsegona
                      \advance\fletxadosdos by \difsegona
                 \else\advance\fletxauu by -\difsegona
                      \advance\fletxaudos by -\difsegona
          \fi
    \diftercera=\wd\Sis
    \advance\diftercera by -\wd\Nou
    \divide\diftercera by 2
    \ifdim\diftercera>0pt\advance\fletxadosdos by \diftercera
                  \else\advance\fletxaudos by -\diftercera
          \fi
    $$
    \align
    \ifdim\wd\Una>0pt\copy\Una \colon \,\,
              \else\fi
    & \copy\Dues\vbox{\hsize=\fletxauu
                     \offinterlineskip
                     \hbox to \fletxauu{\hfil\copy\Tres\hfil}
                     \kern0pt
                     \hbox to \fletxauu{\kern1.5pt
                                        \rightarrowfill
                                        \kern1.5pt}}
      \copy\Quatre\vbox{\hsize=\fletxaudos
                      \offinterlineskip
                      \hbox to \fletxaudos{\hfil\copy\Cinc\hfil}
                      \kern0pt
                      \hbox to \fletxaudos{\kern1.5pt
                                           \rightarrowfill
                                           \kern1.5pt}}
      \copy\Sis\\
    & \hskip\difprimera
      \copy\Set\hbox to \fletxadosu{\kern1.5pt
                                    \mapstofill
                                    \kern1.5pt}
      \copy\Vuit\hbox to \fletxadosdos{\kern1.5pt
                                      \mapstofill
                                      \kern1.5pt}
      \copy\Nou
    \endalign
    $$}
\def\aplicaciotag: nom: #1 conjunt1: #2 aplicacio: #3
     conjunt2: #4 element1: #5
     element2: #6 tag: #7 fi{
    \fletxauu=1.2cm
    \fletxadosu=1.2cm
    \setbox\Una=\hbox{#1}
    \setbox\Dues=\hbox{#2}
    \setbox\Tres=\hbox{#3}
    \setbox\Quatre=\hbox{#4}
    \setbox\Cinc=\hbox{#5}
    \setbox\Sis=\hbox{#6}
    \difprimera=\wd\Tres
    \advance\difprimera by-0.4cm
    \ifdim\difprimera>0pt\advance\fletxauu by \difprimera
                         \advance\fletxadosu by \difprimera
                         \fi
    \difprimera=\wd\Dues
    \advance\difprimera by -\wd\Cinc
    \divide\difprimera by 2
    \ifdim\difprimera>0pt\advance\fletxadosu by \difprimera
                    \else\advance\fletxauu by -\difprimera
          \fi
    \difsegona=\wd\Quatre
    \advance\difsegona by -\wd\Sis
    \divide\difsegona by 2
    \ifdim\difsegona>0pt\advance\fletxadosu by \difsegona
                   \else\advance\fletxauu by -\difsegona
          \fi
    $$
    \eqalign{
    \ifdim\wd\Una>0pt\copy\Una \colon \,\,
              \else\fi
     & \copy\Dues \vbox{\hsize=\fletxauu
                             \offinterlineskip
                             \hbox to \fletxauu{\hfil \copy\Tres \hfil}
                             \kern0pt
                             \hbox to \fletxauu{\kern1.5pt
                                                \rightarrowfill
                                                \kern1.5pt}}
                \copy\Quatre \cr
              & \hskip\difprimera
                \copy\Cinc \hbox to \fletxadosu{\kern1.5pt
                                            \mapstofill
                                            \kern1.5pt}
                \copy\Sis
    }
    \eqno (#7)
    $$}
\def\aplicacio: nom: #1 conjunt1: #2 aplicacio: #3
     conjunt2: #4 element1: #5
     element2: #6 fi{
    \fletxauu=1.2cm
    \fletxadosu=1.2cm
    \setbox\Una=\hbox{#1}
    \setbox\Dues=\hbox{#2}
    \setbox\Tres=\hbox{#3}
    \setbox\Quatre=\hbox{#4}
    \setbox\Cinc=\hbox{#5}
    \setbox\Sis=\hbox{#6}
    \difprimera=\wd\Tres
    \advance\difprimera by-0.4cm
    \ifdim\difprimera>0pt\advance\fletxauu by \difprimera
                         \advance\fletxadosu by \difprimera
                         \fi
    \difprimera=\wd\Dues
    \advance\difprimera by -\wd\Cinc
    \divide\difprimera by 2
    \ifdim\difprimera>0pt\advance\fletxadosu by \difprimera
                    \else\advance\fletxauu by -\difprimera
          \fi
    \difsegona=\wd\Quatre
    \advance\difsegona by -\wd\Sis
    \divide\difsegona by 2
    \ifdim\difsegona>0pt\advance\fletxadosu by \difsegona
                   \else\advance\fletxauu by -\difsegona
          \fi

    $$
    \begin{aligned}
    \ifdim\wd\Una>0pt\copy\Una \colon \,\,
              \else\fi
     & \copy\Dues \vbox{\hsize=\fletxauu
                             \offinterlineskip
                             \hbox to \fletxauu{\hfil \copy\Tres \hfil}
                             \kern0pt
                             \hbox to \fletxauu{\kern1.5pt
                                                \rightarrowfill
                                                \kern1.5pt}}
                \copy\Quatre \\
              & \hskip\difprimera
                \copy\Cinc \hbox to \fletxadosu{\kern1.5pt
                                            \mapstofill
                                            \kern1.5pt}
                \copy\Sis
    \end{aligned}
    $$}

\begin{document}
\title{The Modelling of Random Phenomena}
\author{Aureli Alabert \\ 
           Department of Mathematics \\
           Universitat Autònoma de Barcelona \\  
           08193 Bellaterra, Catalonia \\  
           \url{Aureli.Alabert@uab.cat}  
}

\thispagestyle{empty}
\maketitle

\tableofcontents

\section{Introduction to randomness}

\subsection{Random phenomena}
    A \Conc{random phenomenon} is a physical phenomenon in which 
    ``randomness'' takes a place.
    
    So, what is \Conc{randomness}?    
    It is something that we do not control, in the sense that it may lead to different
    outcomes or measurements of the phenomenon in what we believe are ``identical'' 
    conditions.
    
    There are many keywords associated to the discussion and mathematical foundation of random phenomena:
    probability, chance, likelihood, statistical regularity, 
    plausibility, \dots 
    There are whole books discussing and trying to explain what is the nature of chance and randomness.
    It is not worth going into such philosophical depth for the practitioner. One may get lost 
    into the variety of ''definitions'' or ``trends'' related to the word \Conc{probability}
    (classical, frequentist, axiomatic, subjective, objective, logical, \dots) 
    or \Conc{statistics} (frequentist, classical, Bayesian, decision-theoretic, \dots).    

\subsection{The modelling point of view}
  Instead, take the modelling point of view: Each problem
  must be treated in its own merits, choosing the appropriate tools provided by mathematics.
    
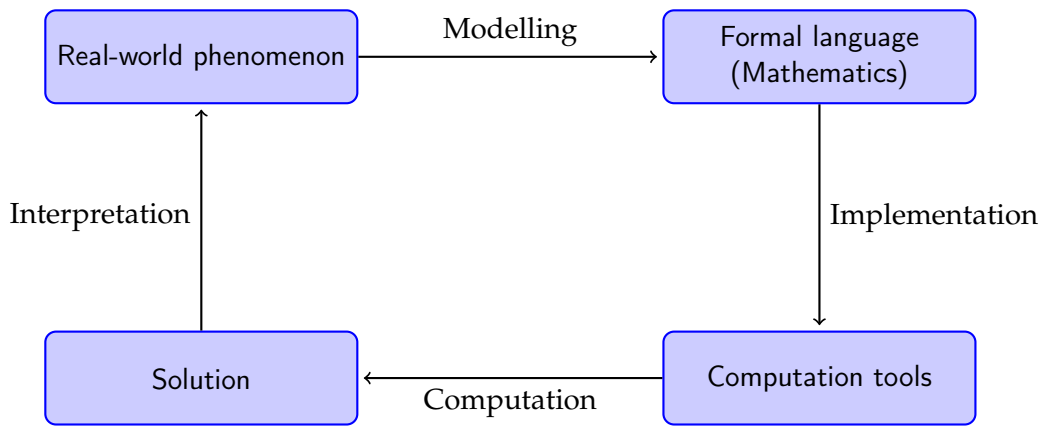
\begin{figure}
  \begin{center}
    \usetikzlibrary{shapes.geometric}
    \usetikzlibrary{shapes.arrows}
    \begin{tikzpicture} [
    auto,
    block/.style    = { rectangle, draw=blue, thick, 
      fill=blue!20, text width=10em, text centered,
      rounded corners, minimum height=3.2em },
    line/.style     = { draw, thick, ->, shorten >=2pt },
    ]
    \matrix [column sep=40mm, row sep=30mm] {
      \node [block] (rworld) {\textsf{Real-world phenomenon}};   &  
      \node [block] (maths) {\textsf{Formal language \\ (Mathematics)}};\\
      \node [block] (sol) {\textsf{Solution}}; &
      \node [block] (comp) {\textsf{Computation tools}}; \\
    };
    \begin{scope} [every path/.style=line]
    \path (sol)      --    node [midway] {Interpretation} (rworld);
    \path (rworld)      --   node [midway] {Modelling} (maths);
    \path (maths)   -- node [midway] {Implementation} (comp);
    \path (comp)   --    node [midway] {Computation} (sol);
    \end{scope}
    \end{tikzpicture}
  \end{center}  \caption{General mathematical modelling}\label{fig:mod_scheme}
\end{figure}

  In general, the modelling of a real world phenomenon follows the scheme of Figure 
  \ref{fig:mod_scheme}.  

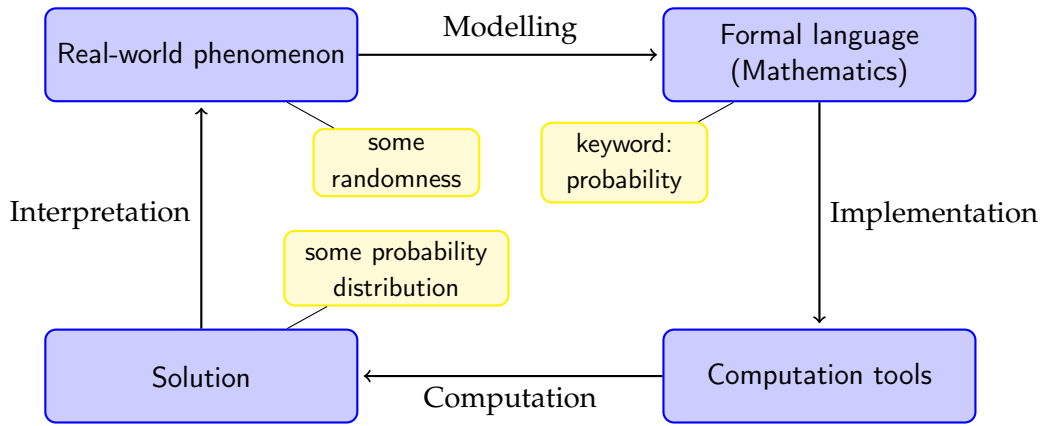
\begin{figure}
  \begin{center}
    \usetikzlibrary{shapes.geometric}
    \usetikzlibrary{shapes.arrows}
    \begin{tikzpicture} [
    auto,
    block/.style    = { rectangle, draw=blue, thick, 
      fill=blue!20, text width=10em, text centered,
      rounded corners, minimum height=3.2em },
    block2/.style    = { rectangle, draw=yellow, thick, 
      fill=yellow!20, text width=5em, text centered,
      rounded corners, minimum height=1.5em },
    block3/.style    = { rectangle, draw=yellow, thick, 
      fill=yellow!20, text width=7em, text centered,
      rounded corners, minimum height=1.5em },
    line/.style     = { draw, thick, ->, shorten >=2pt },
    segment/.style ={ draw, -},
    ]
    \matrix [column sep=40mm, row sep=30mm] {
      \node [block] (rworld) {\textsf{Real-world phenomenon}};   &  
      \node [block] (maths) {\textsf{Formal language \\ (Mathematics)}};\\
      \node [block] (sol) {\textsf{Solution}}; &
      \node [block] (comp) {\textsf{Computation tools}}; \\
    };
    \begin{scope} [every path/.style=line]
    \path (sol)      --    node [midway] {Interpretation} (rworld);
    \path (rworld)      --   node [midway] {Modelling} (maths);
    \path (maths)   -- node [midway] {Implementation} (comp);
    \path (comp)   --    node [midway] {Computation} (sol);
    \end{scope} 
    \begin{scope} [every path/.style=segment]
    \node [block2] (rnd) at (-1.5, 0.7){\small\textsf{some\\randomness}};
    \path (rnd) -- (rworld);
    \node [block2] (kword) at (1.5, 0.7){\small\textsf{keyword:\\probability}};
    \path (kword) -- (maths);
    \node [block3] (dist) at (-1.5, -0.7){\small\textsf{some probability\\distribution}};
    \path (dist) -- (sol);
    \end{scope}
    \end{tikzpicture}
  \end{center}  \caption{Mathematical modelling in the presence of randomness}\label{fig:mod_scheme2}
\end{figure}

  When randomness is present, the scheme is the same. The distinguishing feature is the 
  use of the mathematical concept of ``probability'' (which has an unambiguous and worldwide
  accepted definition), and the solution to
  the problem comes usually in the form of a ``probability distribution'' or some particular
  property of a probability distribution. See Figure \ref{fig:mod_scheme2}.

\subsection{Quantifying randomness: Probability}
   
\begin{figure}
\begin{center}
   \begin{overpic}
     [scale=0.40,tics=10]{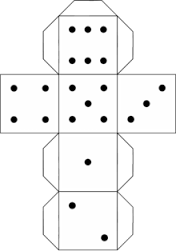}
   \end{overpic}
\end{center}  \caption{A playing die developed to show all its faces.}
\label{fig:die}
\end{figure}

  Take a playing die, for example (Figure \ref{fig:die}). Throwing a die is a familiar random phenomenon. 
  We need the outcome to be unpredictable (thus potentially different) each time;
  otherwise the die is not useful for playing. On the other hand,  
  the experiment is performed each time in identical conditions: We throw the die on the table
  so that it rebounds several times before stopping. Of course, the conditions 
  are no ``truly'' 
  identical; in this case, our ignorance about the exact physical conditions provides
  the desired unpredictability, therefore the randomness.

Suppose we examine the die, and we see that it looks new, homogeneous, balanced
and with no visible manufacturing defect. Is there any outcome that looks more
likely to appear than some other? If not, then it is logical that any attempt to 
quantify the likelihood of the outcomes lead to assign the same quantity to all outcomes.

\begin{figure}
  \begin{center}                                                    
\begin{tikzpicture}
\filldraw[color=red!60, fill=red!5, very thick](0,0) circle (3);
\draw [thick, ] (0,-3.3) -- (0,3.3);
\draw [thick, ] (0,0) --  (30:3.3);
\draw [thick, ] (0,0) --  (210:3.3);
\draw [thick, ] (0,0) --  (150:3.3);
\draw [thick, ] (0,0) --  (330:3.3);
\draw (60:3.5) node {\Cube{1}};
\draw (0:3.5) node {\Cube{2}};
\draw (300:3.5) node {\Cube{3}};
\draw (240:3.5) node {\Cube{4}};
\draw (180:3.5) node {\Cube{5}};
\draw (120:3.5) node {\Cube{6}};
\draw (60:2) node {{\color{blue}$\frac{1}{6}$}};
\draw (0:2) node {{\color{blue}$\frac{1}{6}$}};
\draw (300:2) node {{\color{blue}$\frac{1}{6}$}};
\draw (240:2) node {{\color{blue}$\frac{1}{6}$}};
\draw (180:2) node {{\color{blue}$\frac{1}{6}$}};
\draw (120:2) node {{\color{blue}$\frac{1}{6}$}};
\end{tikzpicture}
\end{center}  \caption{A (presumed) balanced die eating the probability cake.}
\label{fig:diePie}
\end{figure}
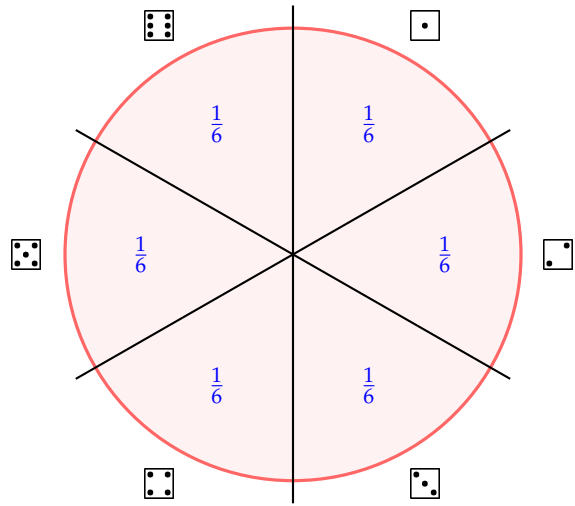

We may think that every outcome takes an equal part of a cake they have to share.
Let us say, arbitrarily, that the cake measures 1. Therefore, every outcome has to
take 1/6 of the cake. We say that every possible \Conc{result} $\omega$ of the 
random phenomenon ``throwing a balanced die'' has a \Conc{probability} of 1/6.
See Figure \ref{fig:diePie}.
  
    From the probability of all possible \Conc{results} $\omega\in \Omega$, we can deduce 
    (define, in fact, but in the only sensible way)
    the probability of
    all possible \Conc{events}, that is, subsets $A\subset\Omega$: The event $A$ takes the part of cake
    that its results $\omega\in A$ take in total.

\subsection{The law of Large Numbers}
    \Conc{The relative frequency} of an event in a series of identical experiments is the quotient 
\begin{equation*}
  \frac{\text{Number of occurrences of the event}}{\text{Number of experiments performed}} 
\end{equation*}

    If $1/6$ is the probability of  
    obtaining a $3$ when tossing the die, \emph{it can be proved} that the relative frequency of the event $\{3\}$
    converges to $1/6$ when the number of experiments tends to infinity.
    
    In general, the relative frequency of an event converges to its probability. This is the 
    \Conc{Law of Large Numbers}. It is a Theorem (an important one). It is not a definition of ``probability'',
    as it is frequently said.

\subsection{Statistical inference}
  We may think that a die is balanced when in fact it is not. In this case, the relative frequencies
  will not converge to the probabilities that we expect. Or,
  plainly, we suspect that the die is not balanced, and we do not know what to expect.

  In any case, the Law of Large Numbers leads to the following idea:     
    \begin{enumerate}
    \item
    Toss the die as many times as you can. 
    \item
    Write down the relative frequency of each result. 
    \item
    Construct the model of the die by assigning 
\begin{equation*}
  \text{Probability of $\omega$ := Relative frequency of $\omega$ .} 
\end{equation*}    
    \end{enumerate}
    
  This is \Conc{Statistical Inference}: We construct a model of a random phenomenon 
  using the data provided by a \Conc{sample} of the \Conc{population}. 
  
  The population here is a (hypothetical) infinite sequence of die throws. In the usual applications,
  the population is a big, but finite, set of objects (people, animals, machines or anything), and the
  sample is a subset of this set.

  In another common (and definitely overused) setting of statistical inference, 
  one simply declares the
  die as balanced unless the relative frequencies deviate too much of the expected
  values. If they do, then the die is declared ``non-balanced''.

\subsection{Probability. The mathematical concept}
 
    We want a mapping that assigns to every event a number called ``the probability of the event''
    satisfying:
    
    \begin{enumerate}
    \item  
    It is nonnegative. 
    \item  
    The probability of the whole set  $\Omega$ of possible results is 1. 
    \item
    The probability of the union of two disjoint events is the sum of the 
    probabilities of the two events.
    \end{enumerate}
    
      Formally:
    A \Conc{probability} is a mapping  
\aplicacio:
  nom: $P$
  conjunt1: ${\cal P}(\Omega)$
  aplicacio: {}
  conjunt2: $[0,1]$
  element1: $A$
  element2: $P(A)$
  fi
  
  such that $P(\Omega)=1$ and for any countable family $\{A_n\}_n \subset \Omega$, with $A_i\cap A_j =\emptyset$ if $i\not= j$,
\begin{equation*}
  P\Big(\unio_{n=1}^{\infty} A_n\Big) = \sum_{n=1}^{\infty} P(A_n)
  \ .
\end{equation*}  
  This definition captures perfectly the idea of the pieces of cake taken by the different events that we saw in 
  Figure \ref{fig:diePie}. The extension to a countably infinite union instead of just finite does not
  harm and allows to construct a mathematical theory much more in line with the phenomena that we intend
  to model. Demanding the same for uncountable unions, on the contrary, would collapse the theory and make it useless. 
  If $\Omega$ is
  a finite set, then of course this discussion is void.

    Sometimes it is not possible to define the mapping on the whole collection ${\cal P}(\Omega)$ 
    of subsets of $\Omega$ preserving at the same time 
    the properties of the definition.
    In this case, we define it on a subcollection ${\cal F}\subset {\cal P}(\Omega)$ satisfying some
    desirable stability properties:
    
  \begin{enumerate}
  \item 
    $\Omega \in {\cal F}$ 
  \item
    $A \in {\cal F}$ $\Rightarrow$ $A^c \in {\cal F}$ 
  \item
    $\{A_n\}_n \subset {\cal F}$ $\Rightarrow$ 
    $\unio_{n=1}^{\infty} A_n \in {\cal F}$ ,
  \end{enumerate}  
  where $A^c:=\Omega-A$ is the complement set of $A$.
  
  These subcollections are called \Conc{$\sigma$-fields} or \Conc{$\sigma$-algebras}. 
  They enjoy the right stability properties
  so that the additivity property in the definition of $P$ still makes sense.
  
  \Conc{Probability Theory} is a specialised part of 
  \Conc{Measure and Integration Theory}. In general, a \Conc{measure}
  is a function defined on the sets of a $\sigma$-field with values in a set which is not necessarily
  the interval $[0,1]$.

\subsection{Drawing probabilities}
\begin{figure}
  \begin{center} 
  \begin{tikzpicture}[
  filled/.style={fill=blue!20, draw=blue!50, thick},
  outline/.style={draw=blue!50, thick}]
  \begin{scope}
  \clip (0,0) circle (1.7cm);
  \fill[filled] (0:2cm) circle (1.5cm);
  \end{scope}
  \draw[outline] (0,0) circle (1.7cm) node {$A$};
  \draw[outline]  (0:2cm) circle (1.5cm) node {$B$};
  \draw[outline, black, thick] (-3,-3) rectangle (5,3);
  \node at (5.3,3) {$\Omega$};
  \end{tikzpicture}
\end{center}
\caption{Probabilities and areas}
\label{fig:regions}
\end{figure}
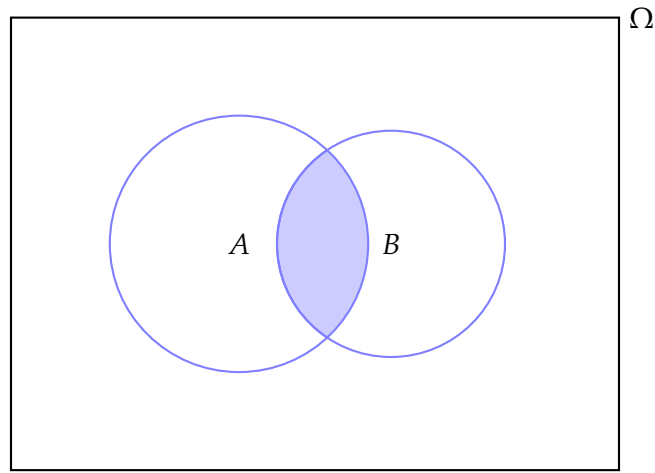

     Probabilities behave like areas of planar regions. Consider Figure \ref{fig:regions}.

  To compute the area of the region $A\cup B$, we may add the areas of $A$ and $B$, and then
  subtract the area of $A\cap B$, which have been counted twice. This leads immediately to the fundamental
  formula:
\begin{align*}
P(A\cup B) = P(A)+P(B)-P(A\cap B)
\ .
\end{align*}
  All usual lists of ``properties of the probabilities'' are trivial derivations of this formula, and
  can also be deduced from Figure \ref{fig:regions}. It is useless to learn them by heart.

\subsection{Conditional probabilities}\label{sec:condprob}  
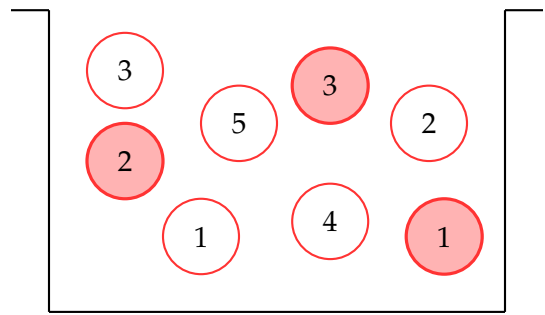
\begin{figure}[t]
\begin{center}  
  \begin{tikzpicture}[
  filled/.style={fill=red!30, draw=red!80, thick},
  outline/.style={draw=red!80, thick}]
  \draw[outline] (2,1) circle (0.50cm);
  \node at (2,1) {$1$};
  \draw[outline] (1,3.2) circle (0.50cm);
  \node at (1,3.2) {$3$};
  \draw[outline] (3.7,1.2) circle (0.50cm);
  \node at (3.7,1.2) {$4$};
  \draw[outline] (2.5,2.5) circle (0.50cm);
  \node at (2.5,2.5) {$5$};
  \draw[outline] (5,2.5) circle (0.50cm);
  \node at (5,2.5) {$2$};
  \filldraw[color=red!80, fill=red!30, very thick](1,2) circle (0.5);
  \draw (1,2) node {$2$};
  \filldraw[color=red!80, fill=red!30, very thick](3.7,3) circle (0.5);
  \draw (3.7,3) node {$3$};
  \filldraw[color=red!80, fill=red!30, very thick](5.2,1) circle (0.5);
  \draw (5.2,1) node {$1$};
  \draw[thick] (0,0) -- (6,0);
  \draw[thick] (6,0) -- (6,4);
  \draw[thick] (0,0) -- (0,4);
  \draw[thick] (0,4) -- (-0.5,4);
  \draw[thick] (6,4) -- (6.5,4);
  \end{tikzpicture}
\end{center}
\caption{Balls in a box}
\label{fig:urn}
\end{figure}

  Consider the following example (see Figure \ref{fig:urn}): We have a box with five white balls, numbered
  1 to 5, and three red balls, numbered 1 to 3.   We pick a ball ``completely at random''. 
  What is the probability of drawing an even number? 

  First of all, what is the probability of picking a particular ball? The expression ``completely at random'',
  though imprecise, is used to mean that all outcomes are equally likely, as
  in the case of the balanced die. 
  
  We are interested in the event
  $A=\{W_2, W_4, R_2\}$, where $W$ means white ball and $R$ red ball. Since each of the balls in $A$ takes
  $\frac{1}{8}$ of the probability cake, we have that
\begin{equation*}
  P(A)=\frac{1}{8}+\frac{1}{8}+\frac{1}{8} = \frac{3}{8}
  \ .
\end{equation*}
  Now suppose a ball has been picked by someone, who tell us that the ball is white. 
  What is the probability that the ball carries an even number?

  In this case the possible results are $W=\{W_1,W_2,W_3,W_4,W_5\}$, all with probability
    $\frac{1}{5}$, thus the probability of $\{W_2,W_4\}$ is $\frac{2}{5}$. The additional information 
    has led as to change the model, and consequently the value of the probabilities. 
      
Notice that : 
    \begin{align*}
    \frac{2}{5} =
    \frac{2/8}{5/8} =    
    \frac{P(A\cap W)}{P(W)}
    \ ,
    \end{align*}
   where the probabilities in the quotient are those of the original model.    

    The \Conc{conditional probability} of $A$ to $B$ is defined as 
    \begin{align*}
    P\big(\condprob{A}{B}\big) :=  \frac{P(A\cap B)}{P(B)}
    \ .
    \end{align*}
    
    In relation to Figure \ref{fig:regions}, the conditional probability of $A$ to $B$ is the
    proportion of area of $A$ inside $B$.  
    
    We say that $A$ and $B$ are \Conc{independent} if the information that $B$ has happened does
    not change the probability of $A$:
    \begin{equation*}
    P\big(\condprob{A}{B}\big)=P(A)
    \ .
    \end{equation*}
    Equivalently,
    \begin{equation*}
    P(A\cap B)=P(A)\cdot P(B)
    \ .
    \end{equation*}

\subsection{Random variables}  
  We can now step into a second level of difficulty: the concept of random variable. 
  Let us consider the following
  example: We toss two balanced dice, and we are interested in the sum of the points shown. 
  We may consider
  directly the set $\Omega=\{2,\dots,12\}$ and assign probabilities to each element of $\Omega$, 
  but this is difficult; 
  or we may keep the model closer to the real experiment by defining  
  $\Omega=\{(i,j):\ 1\le i\le 6,\ 1\le j\le 6\}$,    
  and think of the mapping
    \aplicacio:
  nom: {}
  conjunt1: $\Omega$
  aplicacio: {$X$}
  conjunt2: $\{2,\dots,12\}$
  element1: $(i,j)$
  element2: $i+j$
  fi
  
    If the dice really look balanced, and if it is clear that the outcome of one die does not influence the outcome
    of the other, then it is natural to distribute the same amount of the probability cake to 
    every pair $(i,j)$, that means $P\{(i,j)\}=\frac{1}{36}$.
    
    This setting induces a probability $P_X$ on $\{2,\dots,12\}$, which is what we are looking for: 
    \begin{align*}
    &P_X\{2\}=P\{(1,1)\}=\frac{1}{36}
    \\
    &P_X\{3\}=P\{(1,2),(2,1)\}=\frac{2}{36}
    \\
    &P_X\{4\}=P\{(1,3),(2,2),(3,1)\}=\frac{3}{36}  \quad \dots\ ,\ \text{etc}
    \end{align*}

  In general, 
    a \Conc{random variable} is a mapping $X\colon \Omega \longrightarrow \Reals$.
    ($\Reals$ can be replaced by other convenient sets; technically, the random variable must take 
    values in another \Conc{measure space}, that is, a set endowed with a $\sigma$-field.) 
    The \Conc{law} of a random variable is the 
    probability $P_X$ on $\Reals$ induced by
    $P$ and $X$ as in the example.

    From the modelling point of view \emph{the law is the important thing}, not $\Omega$ or the mapping $X$
    themselves. Typically one says: ``I am observing a random phenomenon following the law \dots''.

    From the law of a random variable one may define certain numeric values that carry some 
    information, and that sometimes are all that is needed in a particular application. The most 
    important one is the \emph{expectation}, which is the ``mean value'' that a variable with 
    that law will take. It can be though as the limit of the arithmetic mean of the observed
    values of the variable when the number of observations tends to infinity. But this is 
    again a version of the Law of Large Numbers, and not a definition. 
    
    The \Conc{expectation} 
    $\E[X]$ of a random variable $X$ with law $P_X$ is defined as 
\begin{equation*}
    \E[X]:=\sum k \cdot P_X\{k\}
    \ ,
\end{equation*}
   with the sum extended over all values taken by $X$.  
    The \Conc{variance} of $X$ is a degree of dispersion of its values around the expectation,
    and defined as
\begin{equation*}
   \Var[X]:=\E\big[(X-\E[X])^2\big]
   \ .
\end{equation*}     
\subsection{The binomial law}\label{binomiallaw} 
  Leaving aside the elementary ``equiprobable'' or ``uniform'' model of the balanced die, the most basic useful
  example of probability law is the one appearing in the following situation:

  Fix an event $A$ of any random experiment. Call $p$ its probability: $P(A)=p$. 
  Repeat $n$ times the same 
  experiment, and let 
  $X$ be the number of occurrences of $A$ in the $n$ trials. The law of $X$ is then determined
  by 
\begin{equation}\label{eq:binomiallaw}
   P\{X=k\} = 
   {n \choose k} p^k (1-p)^{n-k}
   \ ,\quad
   k=0,\dots,n
   \ .
\end{equation}
    We write $X\sim \Binom(n,p)$ and say that $X$ follows a \Conc{binomial law} with parameters 
    $(n,p)$.
    
   The sentence ``repeating $n$ times the \emph{same} experiment'' means in 
   particular that one experiment may not influence the result of another, 
   and therefore events concerning the outcome or one experiment are independent
   of events concerning the outcome of the other experiments, in the sense 
   of section  \ref{sec:condprob}. This fact is key in the deduction of 
   formula (\ref{eq:binomiallaw}).

\section{Examples from daily life: Arrivals and waiting lines}

\subsection{The geometric law}  
    Assume the experiments of Section \ref{binomiallaw} are performed continuously and 
    at regular unit time intervals. 
    We want to know 
    the time elapsed between an occurrence of $A$ and the next occurrence of $A$.  
    Or, in other words, how many experiments are needed before observing again 
    the event $A$.
    
    This is a situation that may be of interest in manufacturing, where the event
    $A$ is the occurrence of a defective item in the production line.
    
    Let $N$ be the number of $A^c$ occurrences before the next occurrence of $A$.     
    Then it is easy to deduce    
  \begin{equation*}
     P\{N=k\} = 
     (1-p)^{k}\cdot p 
     \ ,\quad
     k=0,1,2,\dots
  \end{equation*}
   
    We write $N\sim \Geom(p)$ and say that $N$ follows a \Conc{geometric law} with 
    parameter $p$.

\subsection{Tails and the memoryless property}  
    Once we know the \Conc{density function} (or \Conc{probability function}) 
    $k\mapsto P\{N=k\}$, 
    we can compute, as in the case of the die, the probability
    of any event $P\{N\in B\}$, where $B$ is any subset of $\Naturals$.   
    In particular, we can compute the right and left \Conc{tails} of the law:
\begin{equation*}
  P\{N> k\}= (1-p)^{k+1} 
  \ ,\quad 
  P\{N\le k\}= 1-(1-p)^{k+1}
  \ .
\end{equation*}
   Because of the (hypothesized) independence between the experiments, 
   the law of $N$ is the same if we define
   $N$ as  the number of $A^c$ occurrences before the first occurrence of $A$.  
   From this fact one can prove the \Conc{memoryless property}:
    \begin{equation*}
    P\big\{\condprob{N>m+k}{N>m}\big\} = P\{N>k\}
    \ .
    \end{equation*}
  In words, knowing that 
 the event has not appeared in the first $k$ experiments, it is not more or less
 likely to appear than if we just start now the sequence.

\subsection{Arrivals at random times: The Poisson law}\label{sec:arrivalsrt}
   Assume now that the arrivals occur at random times instead of regularly.     
    For example, the arrival of customers to a waiting line may correspond to 
    this situation. To be
    precise, assume:
    \begin{enumerate}
    \item
      People arrive alone (never in groups).
    \item
      The probability $p$ that an arrival occurs during a time interval of length $h$ (small) 
      is proportional to $h$:
     \begin{equation*}
     p = \lambda \cdot h
     \end{equation*} 
    \item
      The number of arrivals on disjoint time intervals are independent random variables.
    \end{enumerate}
    
    We would like to know, for instance, the law of the number of arrivals $N_t$ in the interval $[0,t]$, or the 
    number of arrivals per unit time. The hypotheses above are quite suitable for a situation where
    the arrivals can be considered ``completely at random''.
    
    Of course, hypothesis 2 can only hold true in an infinitesimal sense. 
    Strictly speaking, one should say
    $\lim_{h\to 0} p/h = \lambda$.

    Now, divide $[0,t]$ in intervals of length $h=t/n$. 
    For $n$ big enough, inside each interval we will see at most one arrival, 
    and this will happen with probability $\lambda h$. 
    Therefore, the number of arrivals in $[0,t]$ follows
    approximately a law $\Binom(n,\lambda t/n)$. Hence, by (\ref{eq:binomiallaw}):
    
  \begin{align*}
    P\{ \text{$k$ arrivals in $[0,t]$}\}  
    = {n\choose k} \cdot \Big(\frac{\lambda t}{n}\Big)^k \cdot 
    \Big(1-\frac{\lambda t}{n}\Big)^{n-k}
    \ .
  \end{align*}
    Taking $n\to\infty$,
  \begin{align}\label{eq:lawNt}
    P\{ \text{$k$ arrivals in $[0,t]$}\}= 
    \frac{(\lambda t)^k}{k!} \exp\{-\lambda t\}
    \ .
  \end{align}
  
   Let $N$ be the number of arrivals per unit time. We write $N\sim \Pois(\lambda)$ and say 
   that $N$ follows a \Conc{Poisson law} with parameter $\lambda$:
  \begin{align*}
    P\{ N=k\} =
    \frac{\lambda^k}{k!} \exp\{-\lambda\}
    \ .
  \end{align*}
   The parameter $\lambda$ is called the \Conc{traffic intensity}.

\subsection{Interarrival times: The exponential law}  
  Let $T$ be the time between two arrivals. As in the case of the geometric law, 
  this random variable
  is equal in law to the time when the first arrival takes place.
  The event $\{T>t\}$ means to have ${\text{0 arrivals in $[0,t]$}}$, whose probability
  according to (\ref{eq:lawNt}) is $\exp\{-\lambda t\}$.
  
  We observe that this probability is nonzero for all $t\ge 0$, 
  and that it cannot be expressed as the 
  sum of the probability of elementary events. We say that
  the interarrival times follow a \Conc{continuous law}, in contrast with 
  all laws seen so far, called \Conc{discrete laws}.

   In the case of continuous laws,  
    the \Conc{density} is a function $f\colon \Reals \longrightarrow \Reals^+$ such that 
    $P\{T\in[a,b]\}$ is the area under its graph between $a$ and $b$.
    
    \begin{align*}
      P\{T\in[a,b]\}=\int_a^b f
    \end{align*}

   To compute the density of the interarrival times, we observe that
    \begin{align*}
      \int_0^t f = P\{T\in[0,t]\}= 1-\exp\{-\lambda t\}\ ,
    \end{align*}
   so that    
    \begin{align*}
      f(t)=\lambda\cdot \exp\{-\lambda t\}
    \end{align*}
            
    $T\sim \Exp(\lambda)$ is called the \Conc{exponential law} with parameter $\lambda$.

\subsection{Continuous laws}  
  Continuous laws have some features that contrast with those of discrete laws:    
    \begin{itemize}
    \item
      The law is not determined by the probability of the individual outcomes.    
    \item  
      It is the density that determines the law. (This can be said to be true also
      for discrete laws, but the concept of ``density function'' is different.)    
    \item  
      It is not possible to assign a probability to all subsets of the real line 
      (this is not obvious).       
      But we do not need to! It is possible to assign a probability to all intervals,
      and therefore to the members of the 
      minimal $\sigma$-field containing the intervals, which is far more than what 
      we need from a practical point of view.      
    \item  
      Continuous laws show why we cannot ask a probability to be additive for 
      collections
      of arbitrary cardinality. For example: $1=P\{T\ge 0\}\neq \sum_{t\ge 0}P\{T=t\}=0$.
    \item
      The \Conc{expectation} of a variable with a continuous law cannot be defined with sums.
      It is the integral 
\begin{equation*}
  \E[X]:=\int_{-\infty}^{\infty} x f(x)\,dx
  \ ,
\end{equation*}
      where $f$ is the density. Notice however the analogy with the definition
      for discrete laws. In the context of measure theory, the expectation can be 
      expressed in a unified way for all cases.       
    \end{itemize}

The correct name of these laws is \Conc{absolutely continuous}, for mathematical consistency,
but the adverb is frequently dispensed with. ``Continuous'', strictly speaking, 
simply means that the 
so-called \Conc{distribution function} $F(x):=P\{X\le x\}$, which is always 
non-decreasing and right-continuous, is furthermore continuous; whereas 
``absolutely continuous''
refers to the stronger property that the distribution function is a primitive 
of another function, the density:
$F(x)=\int_{-\infty}^x f$.  

\subsection{Poisson arrivals / Exponential times}  

  Still some remarks about the relation between the Poisson and the exponential laws:
    \begin{enumerate}
    \item
    If the interarrival times are $\Exp(\lambda)$, then the arrivals per unit time are 
    $\Pois(\lambda)$.
    \item 
    This situation is called ``completely random arrivals", in the sense that 
    the arrival times   $0 < t_1 < t_2 < \cdots < t_k < t$ have the law of $k$ independent
      uniformly distributed values in $[0,t]$, after sorting them. 
    \item  
    The exponential laws enjoy the same \emph{memoryless property} as the geometric law,
    \begin{align*}
    P\{\condprob{T> t+s}{T> s}\} = P\{T>t\}
    \ ,
    \end{align*}
    and is the only continuous law with this property.  It is a good model for lifetimes of 
    ``ageless devices''; for instance, the lifetime of 
      an electronic device, or living beings in their middle ages, when the death comes 
      from internal or external accidents (electric shocks, heart strokes, \dots).
    \end{enumerate}

\subsection{The Poisson process} 
    The collection of random variables $\{N_t,\ t\ge 0\}$, counting how many arrivals have occurred
    in the time interval $[0,t]$, form the \Conc{Poisson process}.
    
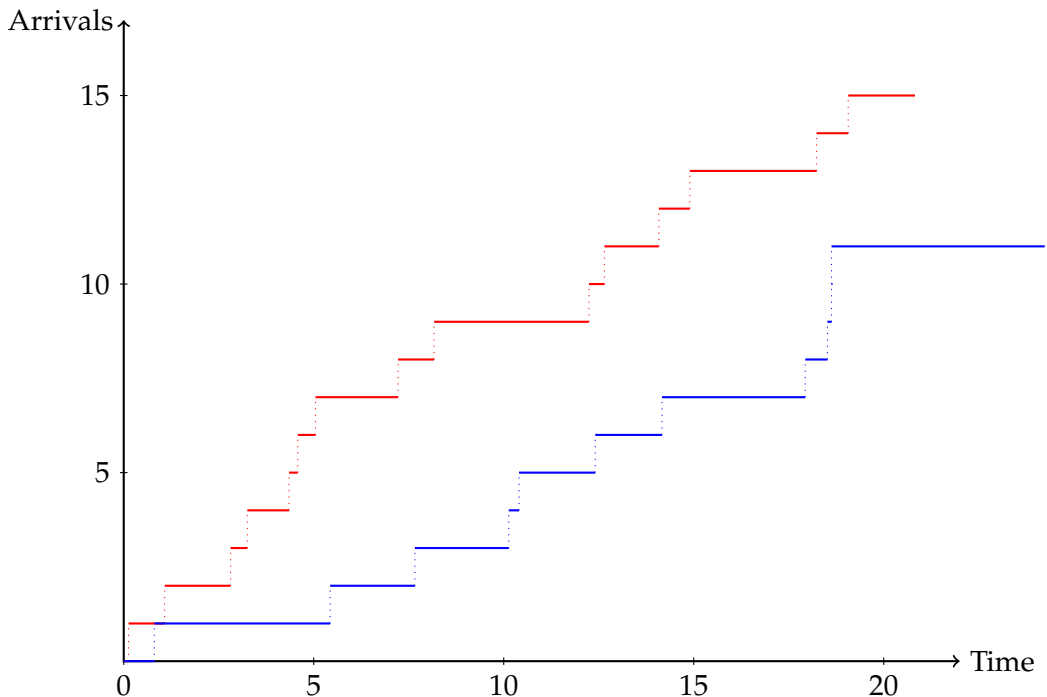
\begin{figure}
  \begin{center}
    \begin{tikzpicture}[scale=0.5] 
    \node [left] at (0,17) {{Arrivals}};
    \node [right] at (22,0) {{Time}};
    \draw [<->, thick] (0,17) -- (0,0) -- (22,0);
    \draw [-] ( -0.1, 5) node [left] {$5$} -- ( 0.1, 5); 
    \draw [-] ( -0.1, 10) node [left] {$10$} -- ( 0.1, 10);  
    \draw [-] ( -0.1, 15) node [left] {$15$} -- ( 0.1, 15); 
    \draw [-] ( 0, -0.1) node [below] {$0$} -- ( 0, 0.1); 
    \draw [-] ( 5, -0.1) node [below] {$5$} -- ( 5, 0.1); 
    \draw [-] ( 10, -0.1) node [below] {$10$} -- ( 10, 0.1); 
    \draw [-] ( 15, -0.1) node [below] {$15$} -- ( 15, 0.1); 
    \draw [-] ( 20, -0.1) node [below] {$20$} -- ( 20, 0.1); 
    
    \draw [-, red, thick] ( 0, 0.0000000)  -- ( 0.1269351, 0);          
    \draw [dotted, red] ( 0.1269351, 0) -- ( 0.1269351, 1);  
    \draw [-, red, thick] (0.1269351, 1)  -- ( 1.0767142, 1);       
    \draw [dotted, red] ( 1.0767142, 1) -- ( 1.0767142, 2);  
    \draw [-, red, thick] ( 1.0767142, 2)  -- ( 2.8113458, 2);      
    \draw [dotted, red] ( 2.8113458, 2) -- ( 2.8113458, 3);  
    \draw [-, red, thick] ( 2.8113458, 3)  -- ( 3.2514947, 3);      
    \draw [dotted, red] ( 3.2514947, 3) -- ( 3.2514947, 4);  
    \draw [-, red, thick] ( 3.2514947, 4)  -- ( 4.3494293, 4); 
    \draw [dotted, red] ( 4.3494293, 4) -- ( 4.3494293, 5);  
    \draw [-, red, thick] ( 4.3494293, 5)  -- ( 4.5802970, 5); 
    \draw [dotted, red] ( 4.5802970, 5) -- ( 4.5802970, 6);  
    \draw [-, red, thick] ( 4.5802970, 6)  -- ( 5.0453761, 6);  
    \draw [dotted, red] ( 5.0453761, 6) -- ( 5.0453761, 7);  
    \draw [-, red, thick] ( 5.0453761, 7)  -- ( 7.2187003, 7);  
    \draw [dotted, red] ( 7.2187003, 7) -- ( 7.2187003, 8);  
    \draw [-, red, thick] ( 7.2187003, 8)  -- ( 8.1683987, 8); 
    \draw [dotted, red] ( 8.1683987, 8) -- ( 8.1683987, 9);  
    \draw [-, red, thick] ( 8.1683987, 9)  -- ( 12.2419776, 9);  
    \draw [dotted, red] ( 12.2419776, 9) -- ( 12.2419776, 10);  
    \draw [-, red, thick] ( 12.2419776, 10)  -- ( 12.6492998, 10); 
    \draw [dotted, red] ( 12.6492998, 10) -- ( 12.6492998, 11);  
    \draw [-, red, thick] ( 12.6492998, 11)  -- ( 14.0817573, 11);  
    \draw [dotted, red] ( 14.0817573, 11) -- ( 14.0817573, 12);  
    \draw [-, red, thick] ( 14.0817573, 12)  -- ( 14.8973963, 12);  
    \draw [dotted, red] ( 14.8973963, 12) -- ( 14.8973963, 13);  
    \draw [-, red, thick] ( 14.8973963, 13)  -- ( 18.2331266, 13); 
    \draw [dotted, red] ( 18.2331266, 13) -- ( 18.2331266, 14);  
    \draw [-, red, thick] ( 18.2331266, 14)  -- ( 19.0650031, 14);  
    \draw [dotted, red] ( 19.0650031, 14) -- ( 19.0650031, 15);  
    \draw [-, red, thick] ( 19.0650031, 15)  -- ( 20.8223631, 15);  

    \draw [-, blue, thick] ( 0, 0.0000000)  -- ( 0.8053696, 0);            
    \draw [dotted, blue] ( 0.8053696, 0) -- ( 0.8053696, 1);  
    \draw [-, blue, thick] ( 0.8053696, 1)  -- ( 5.4323573, 1);       
    \draw [dotted, blue]  ( 5.4323573, 1) -- ( 5.4323573, 2);  
    \draw [-, blue, thick] ( 5.4323573, 2)  -- ( 7.6617429, 2);       
    \draw [dotted, blue] ( 7.6617429, 2) -- ( 7.6617429, 3);  
    \draw [-, blue, thick]  ( 7.6617429, 3)  -- ( 10.1316724, 3);      
    \draw [dotted, blue] ( 10.1316724, 3) -- ( 10.1316724, 4);  
    \draw [-, blue, thick] ( 10.1316724, 4)  -- ( 10.4028206, 4);  
    \draw [dotted, blue] ( 10.4028206, 4) -- ( 10.4028206, 5);  
    \draw [-, blue, thick]  ( 10.4028206, 5)  -- ( 12.4079544, 5);  
    \draw [dotted, blue] ( 12.4079544, 5) -- ( 12.4079544, 6);  
    \draw [-, blue, thick] ( 12.4079544, 6)  -- ( 14.1645107, 6);  
    \draw [dotted, blue] ( 14.1645107, 6) -- ( 14.1645107, 7);  
    \draw [-, blue, thick] ( 14.1645107, 7)  -- ( 17.9344766, 7);  
    \draw [dotted, blue] ( 17.9344766, 7) -- ( 17.9344766, 8);  
    \draw [-, blue, thick] ( 17.9344766, 8)  -- ( 18.5173148, 8);  
    \draw [dotted, blue] ( 18.5173148, 8) -- ( 18.5173148, 9);  
    \draw [-, blue, thick] ( 18.5173148, 9)  -- ( 18.6263642, 9);  
    \draw [dotted, blue] ( 18.6263642, 9) -- ( 18.6263642, 10);  
    \draw [-, blue, thick] ( 18.6263642, 10)  -- ( 18.6263642, 10);  
    \draw [dotted, blue] ( 18.6263642, 10) -- ( 18.6263642, 11);  
    \draw [-, blue, thick] ( 18.6263642, 11)  -- ( 24.2366596, 11);  
    \end{tikzpicture}
  \end{center}
  \caption{A Poisson sample path with $\lambda=1$ (red) and with $\lambda=0.5$ (blue). Lower $\lambda$ means
    less frequent arrivals in average.}\label{fig:poisspath12}
\end{figure}

    When we observe a particular arrival phenomenon, we see, as time passes, 
    a \Conc{sample path} of the Poisson process (see Figure \ref{fig:poisspath12}).
    We may also think of the Poisson process as the collection of all its sample paths.

\subsection{Stochastic processes} 
   In general, a random evolution in time is modelled by a \Conc{stochastic process}. 
   There are two possible points of view of a stochastic process: 
    \begin{enumerate}
    \item
      As a collection of random variables:
     \begin{equation*}
      X:=\{X_t,\ t\ge 0\}\ ,\quad \text{with $X_t\colon \Omega\longrightarrow \Reals$}
      \ .
     \end{equation*}

    \item  
      As a ``random function''
\aplicacio:
  nom: $X$
  conjunt1: $\Omega$
  aplicacio: {}
  conjunt2: $\Reals^{\Reals^+}$
  element1: $\omega$
  element2: $X(\omega)$
  fi
    \end{enumerate}
   Here $\Reals^{\Reals^+}$ denotes the set of all functions $\Reals^+\rightarrow\Reals$, which can 
   be identified  with the Cartesian product of 
   ``$\Reals+$ copies'' of 
   $\Reals$ as a set, as a topological space and as a measure space.

\subsection{Queues (waiting lines)} 
\begin{figure}
\begin{center}
   \begin{overpic}
     [scale=0.5,tics=10]{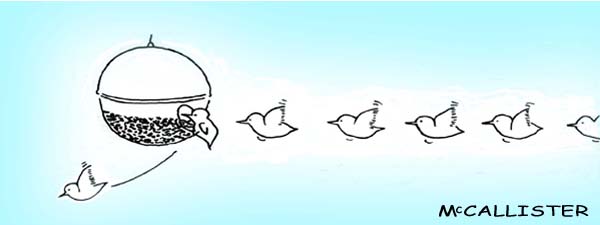} 
   \end{overpic}
\end{center}  \caption{A typical simple queue: Customers arrive, wait in a line, are served, and leave
              the system. (Illustration appeared in The New Yorker, 1977)}\label{fig:bird_queue}
\end{figure}

    A \Conc{queue} is a situation in which \Conc{users} arrive to a \Conc{service}, 
    wait to be served if the service 
    is not immediately available, and leave after having been served (Figure \ref{fig:bird_queue}).
    
    Examples are customers in a supermarket cash, cars in the highway at the toll booths, and parts in a 
    manufacturing chain.
    
    Its behaviour depends, among other things, on:
    \begin{enumerate}
    \item 
      \Conc{Arrival pattern:} Interarrival times, number of users per arrival, 
      patience of the customers, \dots   
    \item  
      \Conc{Service pattern:} Service time, number of users served simultaneously, \dots 
    \item  
      \Conc{Queue discipline:} FIFO (First-In, First-Out), LIFO (Last-In, First-Out), 
      SIRO (Service in Random Order), \dots, with variants specifying 
      priorities, pre-emption, etc.
    \item  
      \Conc{Capacity:} Number of users allowed to wait.
    \end{enumerate}
    
    Moreover, everything may be dependent or not on the state of the system (number of users, etc.) 
    and the elapsed time since the start.

    Typical questions posed in these situations are: 
    \begin{itemize}
    \item
    How many users are in the line? (at a given time, in the mean, \dots)
    \item
    How long a user must wait?  (a given user, in the mean, \dots)
    \item
    How much time a service facility is idle? 
    \item
    How long are the busy/idle  periods of the service facility?
    \end{itemize}

    The answers are random variables if at least one of the features is random. 
    We would like to know the law of these variables, or at least its expectation, or some other
    value of interest.

    The purpose of knowing these laws or law parameters is, frequently, to take a decision about some 
    controllable inputs of the queue, and with some cost associated to each of the values of these
    inputs. For instance, the number of cashiers in a supermarket clearly influences the 
    waiting time of the customers; benefits may increase thanks to that, but the running costs are
    also higher. Here we enter the realm of \Conc{optimisation} and 
    \Conc{operations research}.

\subsection{The M/M/1 queue. Transition probabilities}\label{sec:MM1}
    Assume that we have Poisson arrivals to a queue, the service time is also random and follows an
    exponential law (just one among some common situations), and there is a single service channel 
    (only one user at a time is served).
    
    More precisely, we now put in rigorous mathematics the hypothesis of Section 
    \ref{sec:arrivalsrt}. In the sequel we use the usual notation $o(h)$ to mean any function
    such that $\lim_{h\to 0}o(h)/h=0$.  
    Assume that the  
    arrivals satisfy: 
    \begin{enumerate}
    \item
      $P\{\text{more than one arrival in $[t,t+h]$}\} = o(h)$
    \item  
      $P\{\text{an arrival occurs in $[t,t+h]$}\} = \lambda h + o(h)$
    \item  
      The number of arrivals in non-overlapping time intervals are independent
      random variables.
    \end{enumerate}
    
    And moreover the
      service times satisfy
    \begin{enumerate}
    \item
      $P\{\text{more than one service completed in $[t,t+h]$}\} = o(h)$
    \item  
      $P\{\text{a service is completed in $[t,t+h]$}\} = \mu h + o(h)$  (assuming the service is not idle).
    \item  
      The number of completed services in non-overlapping time intervals are independent
      random variables.
    \end{enumerate}
    
    All these properties together imply that we have a queue where the interarrival times follows the 
    law $\Exp(\lambda)$ and the service times follow the law $\Exp(\mu)$.
    
    Assume, moreover, than jointly considered, arrivals and services are independent.

    Let us call now $N_t$ the number of users in the system at time $t$. We can compute the probability
    that the state of the system changes from $n$ users to any other number in some time interval 
    $[t,t+h]$. These are called the 
    \Conc{transition probabilities}, and can be considered for any stochastic process. 
    It is easy to find, using the hypotheses above
    that for all $n\ge 1$
    \begin{enumerate}[a)]
    \item
    $P\big\{\raise2pt\hbox{$N_{t+h} = n+1$}\ /\ \raise-2pt\hbox{$N_t = n$}\big\}= 
     \lambda h + o(h)$, for $n\ge 0$.
    \item
    $P\big\{\condprob{N_{t+h} = n-1}{N_t = n}\big\}= 
     \mu h + o(h)$, for $n\ge 1$.
    \item
    $P\big\{\condprob{N_{t+h} = n}{N_t = n}\big\}= 
     1-(\lambda+\mu) h + o(h)$, for $n\ge 1$, and \\
      $P\big\{\condprob{N_{t+h} = 0}{N_t = 0}\big\}= 1-\lambda h + o(h)$.
    \item 
    All other transition probabilities are $o(h)$.
    \end{enumerate}

\subsection{The M/M/1 queue. Differential equations}\label{sec:diffeq}  
    Fix two times $s\le t$.
    Denote $p_{nm}(s,t)$ the conditional probability of being in state $m$ at time $t$,
    conditional to be in state $n$ at time $s$. Then, for $m>0$,    
  \begin{align*}
   &p_{nm}(s,t+h) 
   = 
   \sum_{k\in\Naturals} p_{nk}(s,t) \cdot p_{km}(t,t+h) 
   \\
   =&
   p_{nm}(s,t) \cdot p_{mm}(t,t+h)+
   p_{n,m-1}(s,t) \cdot p_{m-1,m}(t,t+h) +
   p_{n,m+1}(s,t) \cdot p_{m+1,m}(t,t+h) 
   + o(h)
   \\
   =&
   p_{nm}(s,t) \cdot \big(1-(\lambda+\mu) h + o(h)\big) +
   p_{n,m-1}(s,t) \cdot \big(\lambda h + o(h)\big) +
   p_{n,m+1}(s,t) \cdot \big(\mu h + o(h)\big) + o(h)
  \end{align*}
  
  Diving by $h$ and taking $h\to 0$, we obtain
  \begin{equation*}
  \frac{d}{dt} p_{nm}(s,t) = -(\lambda+\mu) p_{nm}(s,t) + \lambda p_{n,m-1}(s,t) + \mu p_{n,m+1}(s,t)
  \ .
   \end{equation*}
  Analogously, for $m=0$, one finds
  \begin{equation*}
  \frac{d}{dt} p_{n0}(s,t) = -\lambda p_{n0}(s,t) + \mu p_{n,1}(s,t)
  \ .
  \end{equation*}
    
  This is a countably infinite system of ordinary differential equations for the conditional
  probabilities $p_{nm}(s,t):=P\big\{\condprob{N_{t} = m}{N_s = n}\big\}$, for $s\le t$, and $n,m\in\Naturals$.

  One can also obtain differential equations for the law of $N_t$ itself: 
    Denote $p_{n}(t)=P\{N_t=n\}$.
    For $n>0$,
  \begin{align*}
   \frac{d}{dt} p_{n}(t) 
   &= 
   \frac{d}{dt}\Big(\sum_{k\in\Naturals} p_{k}(0) p_{kn}(0,t)\Big) 
   \\
   &=
   \sum_{k\in\Naturals} p_{k}(0) \Big[  -(\lambda+\mu) p_{kn}(0,t) + \lambda p_{k,n-1}(0,t) + \mu p_{k,n+1}(0,t) \Big]
   \\
   &=
   -(\lambda+\mu) p_{n}(t) + \lambda p_{n-1}(t) + \mu p_{n+1}(t) 
   \ .
  \end{align*}
  
  And, for $n=0$,
  \begin{equation*}
  \frac{d}{dt} p_{0}(t) = -\lambda p_{0}(t) + \mu p_{1}(t)
  \ .
   \end{equation*}

   We get again a countably infinite system of ordinary differential equations.
   The system can be solved exactly but it is difficult and there is a lot of higher 
   mathematics involved.

\subsection{The M/M/1 queue. Steady-state law}   
    In the long run, as $t$ grows, does the law of $N_t$ stabilises?      
    If this is true, then the derivatives in the system of Section 
    \ref{sec:diffeq} must vanish when $t\to\infty$:
    \begin{align*}
      0 &= -\lambda p_{0} + \mu p_{1}
      \\
     0 
     &=
     -(\lambda+\mu) p_{n} + \lambda p_{n-1} + \mu p_{n+1} 
     \ .
    \end{align*}
    
    By induction, 
    \begin{equation*}
    p_n=\Big(\frac{\lambda}{\mu}\Big)^n\cdot p_0
    \end{equation*}
    
    Using as boundary condition $\sum p_n = 1$, we obtain
    \begin{equation*}
    p_0=\frac{1}{\sum_{n=0}^\infty \big(\frac{\lambda}{\mu}\big)^n}
    \ ,
    \end{equation*}
    hence a necessary condition for the existence of a stabilisation is $\lambda < \mu$. 
    Denote $\rho:=\lambda/\mu$. This number is called the \Conc{traffic intensity} of
    the queue.   
    
    If $\rho \ge 1$, no steady-state exists; in fact, the queue tends to grow forever, as more and more users
    accumulate in it. 

    If, on the contrary, $\rho < 1$, then $p_0=1-\rho$, and we get
    \begin{align*}
    p_n=\rho^n(1-\rho)
    \ ,
    \end{align*}
    which is the probability of having $n$ users in the system, in the long run.

  Knowing the law of the number of users in the system in the long run, it is easy to compute:
  \begin{itemize}
  \item  
    The expectation of the number of users $N$ in the system: 
    \begin{equation*}
    \E[N]=\frac{\rho}{1-\rho} \ .
    \end{equation*}
  \item  
    The expectation of the number of customers $N_q$ in the queue: 
    \begin{equation*}
    \frac{\rho^2}{1-\rho} \ .
    \end{equation*}
  \item  
    The law of the waiting time $T_q$ in queue: 
    \begin{align*}
    P\{T_q=0\}&= 1-\rho \ .
    \\
    P\{T_q\le t\}&= 1-\rho \exp\{-\mu(1-\rho)t\} \quad  \text{(for $t>0$).} 
    \end{align*}
  \item
    The expectation of $T_q$: 
    \begin{equation*}
    \E[T_q]=\frac{\lambda}{\mu(\mu-\lambda)} \ .
    \end{equation*}
  \end{itemize}

\subsection{Complex queueing systems. Simulation} 
  The results above are specific of the M/M/1 queue.   
    There are specific results for other types of queues, and there are also some   
    general results. For instance, the relations
    \begin{align*}
    \E[N] &= \lambda \E[T]
    \\
    \E[N_q] &= \lambda \E[T_q]
    \end{align*}
    which one can easily deduce in the M/M/1 queue, 
    are true, no matter the law of arrivals and service times.

However, 
    except for relatively easy queue systems, there is 
    no hope to find analytical results, as computations become intractable very soon. That means that
    in the real world, one can hardly find closed formulae. 
    
    What to do then? One may propose:
  \begin{itemize}
  \item  
    \emph{Idea 1}: Observe the system long enough, take data and do some sort of statistical inference.
  \item   
    \emph{Idea 2}: Simulate the system in a computer, and do statistical inference as well.
  \end{itemize}
   For idea 1 to work, we need the system really running, some mechanism of observation, and a lot
   of time. In practice, we seldom can afford such luxuries.
    For idea 2, on the other hand, we only need, essentially, 
    a mechanism to generate random numbers.
    
   There are very good random number generators embodied in software. Their outcome is not really random,
   but they can fool any detector of ``non-randomness''. Anyway, if the quality of a stream of 
   such pseudo-random numbers is a concern, it is very easy to use a true random number generator
   based in hardware: Nowadays, several internet sites offer potentially infinite streams of true 
   random numbers produced by a quantum device. And such devices are quite cheap, in fact.

\subsection{Birth and death processes} 
 
    A \Conc{birth and death process} $N_t$ takes values in $\Naturals$ and the change across an
    infinitesimal time interval can only be -1, 0, +1: 
    \begin{align*}
       & P\big\{\condprob{N_{t+h} = n+1}{N_t = n}\big\}= \lambda_n\cdot h + o(h)
       \\
       & P\big\{\condprob{N_{t+h} = n-1}{N_t = n}\big\}= \mu_n \cdot h + o(h)
    \end{align*}
  This is a generalisation of the M/M/1 queue model to transition probabilities that may depend on
    the system state.

    The corresponding system of differential equations for the state of the system becomes
  \begin{align*}
   \frac{d}{dt} p_{n}(t) 
   &=
   -(\lambda_n+\mu_n) p_{n}(t) + \lambda_{n-1} p_{n-1}(t) + \mu_{n+1} p_{n+1}(t) 
   \\
  \frac{d}{dt} p_{0}(t) &= -\lambda_0 p_{0}(t) + \mu_1 p_{1}(t)
  \end{align*}
  
  Birth and death processes have been used, for example, to model the varying size of a 
  biological population under given environmental conditions, or to describe the evolution
  of an epidemic.

\section{Example from industry: Inventories}                                  

\begin{figure}[h]
\begin{center}  
   \begin{overpic}
     [scale=0.4,tics=10]{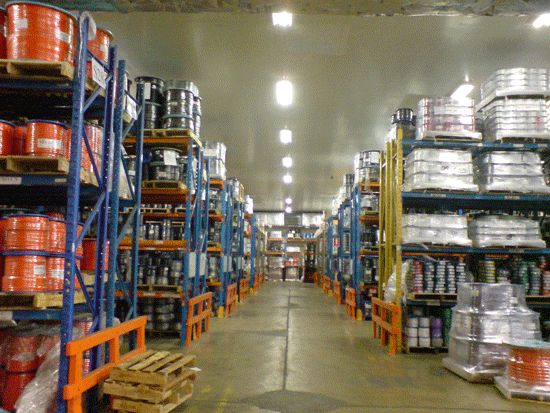}
   \end{overpic}
\end{center}
\caption{A warehouse}
\end{figure}\label{warehouse}

\subsection{Inventory modelling} 
   A company distributes some product, maybe after processing some raw material that arrives
   to the warehouse. Let us assume that we are dealing only with one product and no processing time.
   Assume also that the product has an approximately constant level of demand, but the arrival 
   of orders from the clients is not so predictable.
   The time required to obtain units of product from the manufacturer is also subject to 
   some variability. 
   
   Two fundamental questions in this situation are:  
   
   \begin{enumerate}   
     \item 
     When should more items be ordered?
     \item   
     How many items should be ordered when an order is placed?
   \end{enumerate}
     
   A couple of things to take into account:
   \begin{itemize}
   \item
     If a customer wants to purchase but we do not have items, the sale is lost.
     Therefore, it is important to have enough items in the warehouse.
   \item  
     The product may become obsolete, and there is also a cost of maintaining the inventory.
     Therefore, it is not good to keep in storage too many items.
   \end{itemize}

   Simple hypothesis for an inventory problem that allow analytical computations similar to
   the M/M/1 queue are:
   \begin{itemize}
   \item
     Orders arrive for single items with a random interarrival times following the same law,
     independent from each other.
   \item  
     The time to receive items from the manufacturer (\Conc{lead times}) follows some law,
     and are independent, and independent of order arrival.
   \end{itemize}
          
   A commonly used simple strategy is the \Conc{$(r,s)$-policy}: when the inventory drops to $r$ units, 
   order $s-r$ units. One may measure the performance to this policy, given $r$ and $s$ 
   by the average inventory level, or by the average no-inventory time, or by the number
   of orders that arrive when the inventory is broken, or, most probably,
   by an combination of these and other measures that ultimately reduces to a measure of
   economic benefit that the company wants to maximise.

\begin{figure}
  \begin{center}
    \begin{tikzpicture}[scale=1] 
    \node [left] at (0,6) {{Level}};
    \node [right] at (13,0) {{Time}};
    \draw [<->, thick] (0,6) -- (0,0) -- (13,0);
    \draw [-] ( -0.1, 0) node [left] {$0$} -- ( 0.1, 0); 
    \draw [-] ( -0.1, 1) node [left] {$1$} -- ( 0.1, 1); 
    \draw [-] ( -0.1, 2) node [left] {$2$} -- ( 0.1, 2); 
    \draw [-] ( -0.1, 3) node [left] {$3$} -- ( 0.1, 3); 
    \draw [-] ( -0.1, 4) node [left] {$4$} -- ( 0.1, 4); 
    \draw [-] ( -0.1, 5) node [left] {$5$} -- ( 0.1, 5); 
    \draw [-] ( 0, -0.1) node [below] {$0$} -- ( 0, 0.1); 
    \draw [-] ( 2, -0.1) node [below] {$5$} -- ( 2, 0.1); 
    \draw [-] ( 4, -0.1) node [below] {$10$} -- ( 4, 0.1); 
    \draw [-] ( 6, -0.1) node [below] {$15$} -- ( 6, 0.1); 
    \draw [-] ( 8, -0.1) node [below] {$20$} -- ( 8, 0.1); 
    \draw [-] ( 10, -0.1) node [below] {$25$} -- ( 10, 0.1); 
    \draw [-] ( 12, -0.1) node [below] {$30$} -- ( 12, 0.1); 
    \draw [-, red, thick] ( 0, 5)  -- ( 1, 5);                
    \draw [dotted, red] ( 1,5) -- ( 1,4);  
    \draw [-, red, thick] ( 1, 4)  -- ( 1.6, 4);       
    \draw [dotted, red] ( 1.6,4) -- ( 1.6,3);  
    \draw [-, red, thick] ( 1.6, 3)  -- ( 2.6, 3);      
    \draw [dotted, red] ( 2.6,3) -- ( 2.6,2);  
    \draw [-, red, thick] ( 2.6, 2)  -- ( 3.8, 2);     
    \draw [dotted, red] ( 3.8,2) -- ( 3.8,1);  
    \draw [-, red, thick] ( 3.8, 1)  -- ( 4.2, 1); 
    \draw [dotted, red] ( 4.2,1) -- ( 4.2,4);  
    \draw [-, red, thick] ( 4.2, 4)  -- ( 4.4, 4); 
    \draw [dotted, red] ( 4.4,4) -- ( 4.4,3);  
    \draw [-, red, thick] ( 4.4, 3)  -- ( 4.8, 3); 
    \draw [dotted, red] ( 4.8,3) -- ( 4.8,2);  
    \draw [-, red, thick] ( 4.8, 2)  -- ( 5, 2); 
    \draw [dotted, red] ( 5,2) -- ( 5,5);  
    \draw [-, red, thick] ( 5, 5)  -- ( 6, 5); 
    \draw [dotted, red] ( 6,5) -- ( 6,3);  
    \draw [-, red, thick] ( 6, 3)  -- ( 6.2, 3); 
    \draw [dotted, red] ( 6.2,3) -- ( 6.2,2);  
    \draw [-, red, thick] ( 6.2, 2)  -- ( 6.6, 2); 
    \draw [dotted, red] ( 6.6,2) -- ( 6.6,1);  
    \draw [-, red, thick] ( 6.6, 1)  -- ( 6.8, 1); 
    \draw [dotted, red] ( 6.8,1) -- ( 6.8,0);  
    \draw [-, red, thick] ( 6.8, 0)  -- ( 8, 0); 
    \draw [dotted, red] ( 8,0) -- ( 8,3);  
    \draw [-, red, thick] ( 8, 3)  -- ( 8.2, 3); 
    \draw [dotted, red] ( 8.2,3) -- ( 8.2,2);  
    \draw [-, red, thick] ( 8.2, 2)  -- ( 8.6, 2);  
    \draw [dotted, red] ( 8.6,2) -- ( 8.6,5);  
    \draw [-, red, thick] ( 8.6, 5)  -- ( 8.8, 5);  
    \draw [dotted, red] ( 8.8,5) -- ( 8.8,4);  
    \draw [-, red, thick] ( 8.8, 4)  -- ( 10.2, 4);  
    \draw [dotted, red] ( 10.2,4) -- ( 10.2,3);  
    \draw [-, red, thick] ( 10.2, 3)  -- ( 10.8, 3);  
    \draw [dotted, red] ( 10.8,3) -- ( 10.8,2);  
    \draw [-, red, thick] ( 10.8, 2)  -- ( 11.4, 2); 
    \draw [dotted, red] ( 11.4,2) -- ( 11.4,1);  
    \draw [-, red, thick] ( 11.4, 1)  -- ( 11.8, 1);  
    \draw [dotted, red] ( 11.8,1) -- ( 11.8,0);  
    \draw [-, red, thick] ( 11.8, 0)  -- ( 12.6, 0);  
    \end{tikzpicture}
  \end{center}
  \caption{A path of a inventory process. For some time before day 20 and around day 30 
           the inventory was ``empty''.}\label{inventoryproc}
\end{figure}
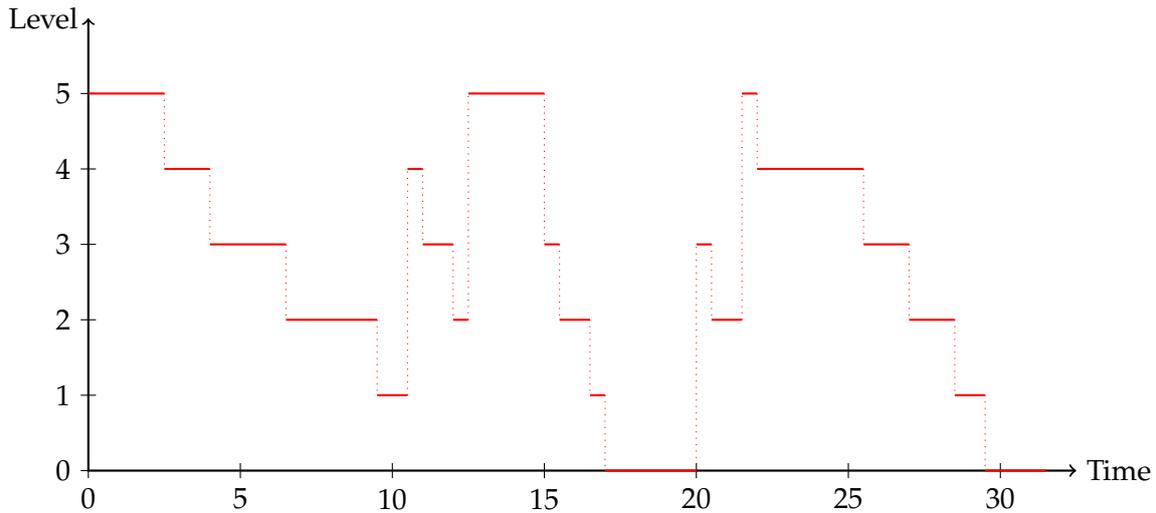

   The inventory process, whose paths have the aspect of Figure \ref{inventoryproc}, 
   is not in general a birth and death process: Items may arrive in batches to the warehouse, and
   the clients' orders may also involve more that one unit.
   It is therefore a generalisation of the situations seen
   in the last sections. But it can still be simulated easily if we know 
   the input distributions.

\subsection{Markov chains} 
  We further generalise by abstracting a remarkable property of the inventory process:
  If we know the state of the system
  at a particular time $t$, we do not need to know anything about previous states
  to predict the future. 
  This is a random analogous of the uniqueness property of deterministic dynamical systems
  when proper initial conditions are given.
  
  Formally: If $t_1<\cdots<t_k<t$, 
   \begin{align*}
       P\big\{\condprob{N_{t} = n}{N_{t_1} = n_1,\ \dots,\ N_{t_k} = n_k}\big\}= 
       P\big\{\condprob{N_{t} = n}{N_{t_k} = n_k}\big\} 
   \end{align*}

   Stochastic processes satisfying this property are called \Conc{Markov chains}, and 
   enjoy an extensive an quite rich theory.

\subsection{Chapman--Kolmogorov equation} 
    Consider times $0\le u<t<s$. Recall the notation of Section \ref{sec:diffeq} for 
    the transition probabilities. 

    The \Conc{Chapman--Kolmogorov equation} for Markov chains establishes that the probability of going
    from state $n$ to state $m$ when time runs from $u$ to $s$ can be computed by
    decomposing all possible paths at the intermediate time $t$:
    \begin{align*}
    \displaystyle p_{nm}(u,s) = \sum_{k} p_{nk}(u,t)p_{km}(t,s)
    \ .
    \end{align*}
  We have already used this in Section \ref{sec:diffeq}.

    In particular, the law of the random variable $N_{t+h}$ can be obtained from 
    the law of $N_t$ and the transition probabilities from $t$ to $t+h$:
    \begin{align*}
    p_{nm}(0,t+h) &= \sum_{k} p_{nk}(0,t)p_{km}(t,t+h)
    \\
    \sum_{n} p_n(0)p_{nm}(0,t+h) &= \sum_{k} \sum_{n} p_n(0)p_{nk}(0,t)p_{km}(t,t+h)
    \\
    p_m(t+h) &= \sum_{k} p_{k}(t)p_{km}(t,t+h)
    \ .
    \end{align*}

\subsection{Kolmogorov forward and backward equations} 
    Assume 
    \begin{align*}
    1-p_{nn}(t,t+h)&=q_n(t)h + o(h)
    \\
    p_{nm}(t,t+h) &= q_{nm}(t)h+o(h) \ , \quad (\text{for $n\not = m$})
    \end{align*}
    for some continuous functions $q_n$ and $q_{nm}$.  
    Then, the following two relations hold:  
    \begin{align*}
    \displaystyle\frac{\partial}{\partial t} p_{nm}(u,t)=q_m(t)p_{nm}(u,t)+ \sum_{k\not=m}p_{nk}(u,t)q_{kj}(t)
    \end{align*}
    \begin{align*}
    \displaystyle\frac{\partial}{\partial u} p_{nm}(u,t)=q_n(u)p_{nm}(u,t)- \sum_{k\not=n}q_{nk}(t)p_{km}(u,t)
    \end{align*}
  These differential equations for the transition probabilities are known as 
  \Conc{Kolmogorov equations}, forward and backward, respectively.

\subsection{Differential equations for the laws} 
    Assume that the functions $q_n$ and $q_{nm}$ above are constant: 
    $q_n(t)\equiv q_n$ and $q_{nm}(t)\equiv q_{nm}$. The Markov chain is then called
    \Conc{time-homogeneous}. 
    
    From Kolmogorov forward equations, letting $u=0$, multiplying by $p_n(0)$ and 
    summing over $n$, one obtains a (infinite) system of differential equations for the laws
    of $N_t$:
    \begin{align*}
    \frac{d}{dt} p_{m}(t)&=-q_m p_{m}(t)+ \sum_{k\not=m}p_{k}(t)q_{kj}
    \end{align*}

\subsection{Long-run behaviour of Markov chains}\label{longrunMC} 
  In many applications it is of interest to study the behaviour of the chain \emph{in the long run}.
  For instance:
\begin{itemize}
  \item
  Limiting distributions: Assume that the limits $\lim_{t\to\infty} p_{nm}(u,t)$ exist and are equal, for all $n$. 
  That means, the limit is independent of the initial state, when time is large. The limit is a probability
  law called the \Conc{limiting} or \Conc{steady-state distribution} of the 
  Markov chain.
  \item
  Stationary distributions: If the limit of the laws $\{\lim_{t\to\infty} p_n(t)\}_n$ exists, it is called the
  \Conc{stationary distribution} of the chain. If there is a limiting distribution, then it coincides with the 
  stationary distribution. But the latter may exist independently. 
  \item
  Ergodicity: Loosely speaking, \Conc{ergodictiy} means that some kind of information that can be extracted from a 
  process as
  a whole, can also be obtained by observing one single path. For instance, ergodicity with respect to the expectation
  means that the limit $\lim_{t\to\infty} \E[X(t)]$ coincides with 
\begin{equation*}
\lim_{t\to\infty} \frac{1}{t}\int_0^t X(s)\,ds 
\end{equation*}  
  for all sample paths $X(s)$. For example, the M/M/1 queue, with traffic intensity $\rho<1$, satisfies this property.
  
  In particular, ergodicity implies that simulating one only sample path for long enough time is sufficient
  to estimate the expectation of the process in the long run.
  \item 
    Classification of states: The elements of the state space of Markov chains are classified according to 
    different interwoven criteria. Among the most important concepts: A state is \Conc{transient} if the 
    probability to never returning to it is positive; otherwise it is called \Conc{recurrent}, and the process
    will certainly visit that state an infinite number of times; a state is \Conc{absorbing} if the chain never
    leaves it once it is reached. 
\end{itemize}      
  
\subsection{Stochastic processes in discrete time} 
  A
  \Conc{discrete time stochastic process} is a process where the family of random 
  variables is indexed by a discrete set, usually $\Enters$ or $\Naturals$.
    
  A \Conc{discrete time Markov chain} has the same
    definition of a Markov (continuous time) chain, except that the index $t$ runs 
    over a discrete set, usually the non-negative integers.

\begin{figure}
\begin{center}
\begin{tikzpicture}[x=1pt,y=1pt]
\definecolor{fillColor}{RGB}{255,255,255}
\path[use as bounding box,fill=fillColor,fill opacity=0.00] (50,50) rectangle (311.35,311.35);
    \node [left] at (50,305) {Passengers};
    \node [right] at (350,50) {Time};
    \draw [<->, thick] (50,300) -- (50,59.83) -- (350,59.83);
    
    {\small
    \draw [-] ( 46, 68.69564) node [left] {$100$} -- ( 50, 68.69564); 
    \draw [-] ( 46, 113.5547) node [left] {$200$} -- ( 50, 113.5547);
    \draw [-] ( 46, 158.4138) node [left] {$300$} -- ( 50, 158.4138);
    \draw [-] ( 46, 203.2729) node [left] {$400$} -- ( 50, 203.2729);
    \draw [-] ( 46, 248.1319) node [left] {$500$} -- ( 50, 248.1319);
    \draw [-] ( 46, 292.991) node [left] {$600$} -- ( 50,  292.991);

    \draw [-] ( 59.83, 55.83) node [below] {$1$} -- ( 59.83, 59.83); 
    \draw [-] ( 82.12, 55.83) node [below] {$13$} -- ( 82.12, 59.83); 
    \draw [-] ( 104.42, 55.83) node [below] {$25$} -- ( 104.42, 59.83); 
    \draw [-] ( 126.72, 55.83) node [below] {$37$} -- ( 126.72, 59.83); 
    \draw [-] ( 149.01, 55.83) node [below] {$49$} -- ( 149.01, 59.83); 
    \draw [-] ( 171.31, 55.83) node [below] {$61$} -- ( 171.31, 59.83); 
    \draw [-] ( 193.60, 55.83) node [below] {$73$} -- ( 193.60, 59.83); 
    \draw [-] ( 215.90, 55.83) node [below] {$85$} -- ( 215.90, 59.83); 
    \draw [-] ( 238.20, 55.83) node [below] {$97$} -- ( 238.20, 59.83); 
    \draw [-] ( 260.49, 55.83) node [below] {$109$} -- ( 260.49, 59.83); 
    \draw [-] ( 282.79, 55.83) node [below] {$121$} -- ( 282.79, 59.83); 
    \draw [-] ( 305.08, 55.83) node [below] {$133$} -- ( 305.08, 59.83); }

\begin{scope}
\path[clip] ( 49.20, 61.20) rectangle (336.15,312.15);
\definecolor{drawColor}{RGB}{255,0,0}

\path[draw=drawColor,line width= 0.4pt,line join=round,line cap=round] ( 59.83, 23.84) -- ( 59.83, 74.08); 

\path[draw=drawColor,line width= 0.4pt,line join=round,line cap=round] ( 61.69, 23.84) -- ( 61.69, 76.77);

\path[draw=drawColor,line width= 0.4pt,line join=round,line cap=round] ( 63.54, 23.84) -- ( 63.54, 83.05);

\path[draw=drawColor,line width= 0.4pt,line join=round,line cap=round] ( 65.40, 23.84) -- ( 65.40, 81.71);

\path[draw=drawColor,line width= 0.4pt,line join=round,line cap=round] ( 67.26, 23.84) -- ( 67.26, 78.12);

\path[draw=drawColor,line width= 0.4pt,line join=round,line cap=round] ( 69.12, 23.84) -- ( 69.12, 84.40);

\path[draw=drawColor,line width= 0.4pt,line join=round,line cap=round] ( 70.98, 23.84) -- ( 70.98, 90.23);

\path[draw=drawColor,line width= 0.4pt,line join=round,line cap=round] ( 72.83, 23.84) -- ( 72.83, 90.23);

\path[draw=drawColor,line width= 0.4pt,line join=round,line cap=round] ( 74.69, 23.84) -- ( 74.69, 84.85);

\path[draw=drawColor,line width= 0.4pt,line join=round,line cap=round] ( 76.55, 23.84) -- ( 76.55, 77.22);

\path[draw=drawColor,line width= 0.4pt,line join=round,line cap=round] ( 78.41, 23.84) -- ( 78.41, 70.49);

\path[draw=drawColor,line width= 0.4pt,line join=round,line cap=round] ( 80.27, 23.84) -- ( 80.27, 76.77);

\path[draw=drawColor,line width= 0.4pt,line join=round,line cap=round] ( 82.12, 23.84) -- ( 82.12, 75.43); 

\path[draw=drawColor,line width= 0.4pt,line join=round,line cap=round] ( 83.98, 23.84) -- ( 83.98, 80.36);

\path[draw=drawColor,line width= 0.4pt,line join=round,line cap=round] ( 85.84, 23.84) -- ( 85.84, 87.09);

\path[draw=drawColor,line width= 0.4pt,line join=round,line cap=round] ( 87.70, 23.84) -- ( 87.70, 84.40);

\path[draw=drawColor,line width= 0.4pt,line join=round,line cap=round] ( 89.56, 23.84) -- ( 89.56, 79.91);

\path[draw=drawColor,line width= 0.4pt,line join=round,line cap=round] ( 91.41, 23.84) -- ( 91.41, 90.68);

\path[draw=drawColor,line width= 0.4pt,line join=round,line cap=round] ( 93.27, 23.84) -- ( 93.27,100.10);

\path[draw=drawColor,line width= 0.4pt,line join=round,line cap=round] ( 95.13, 23.84) -- ( 95.13,100.10);

\path[draw=drawColor,line width= 0.4pt,line join=round,line cap=round] ( 96.99, 23.84) -- ( 96.99, 94.72);

\path[draw=drawColor,line width= 0.4pt,line join=round,line cap=round] ( 98.85, 23.84) -- ( 98.85, 83.50);

\path[draw=drawColor,line width= 0.4pt,line join=round,line cap=round] (100.70, 23.84) -- (100.70, 74.98);

\path[draw=drawColor,line width= 0.4pt,line join=round,line cap=round] (102.56, 23.84) -- (102.56, 86.64);

\path[draw=drawColor,line width= 0.4pt,line join=round,line cap=round] (104.42, 23.84) -- (104.42, 88.89); 

\path[draw=drawColor,line width= 0.4pt,line join=round,line cap=round] (106.28, 23.84) -- (106.28, 91.13);

\path[draw=drawColor,line width= 0.4pt,line join=round,line cap=round] (108.14, 23.84) -- (108.14,103.69);

\path[draw=drawColor,line width= 0.4pt,line join=round,line cap=round] (109.99, 23.84) -- (109.99, 96.96);

\path[draw=drawColor,line width= 0.4pt,line join=round,line cap=round] (111.85, 23.84) -- (111.85,101.00);

\path[draw=drawColor,line width= 0.4pt,line join=round,line cap=round] (113.71, 23.84) -- (113.71,103.69);

\path[draw=drawColor,line width= 0.4pt,line join=round,line cap=round] (115.57, 23.84) -- (115.57,113.11);

\path[draw=drawColor,line width= 0.4pt,line join=round,line cap=round] (117.43, 23.84) -- (117.43,113.11);

\path[draw=drawColor,line width= 0.4pt,line join=round,line cap=round] (119.28, 23.84) -- (119.28,106.38);

\path[draw=drawColor,line width= 0.4pt,line join=round,line cap=round] (121.14, 23.84) -- (121.14, 96.51);

\path[draw=drawColor,line width= 0.4pt,line join=round,line cap=round] (123.00, 23.84) -- (123.00, 89.33);

\path[draw=drawColor,line width= 0.4pt,line join=round,line cap=round] (124.86, 23.84) -- (124.86, 98.31);

\path[draw=drawColor,line width= 0.4pt,line join=round,line cap=round] (126.72, 23.84) -- (126.72,100.55); 

\path[draw=drawColor,line width= 0.4pt,line join=round,line cap=round] (128.57, 23.84) -- (128.57,104.59);

\path[draw=drawColor,line width= 0.4pt,line join=round,line cap=round] (130.43, 23.84) -- (130.43,110.42);

\path[draw=drawColor,line width= 0.4pt,line join=round,line cap=round] (132.29, 23.84) -- (132.29,105.03);

\path[draw=drawColor,line width= 0.4pt,line join=round,line cap=round] (134.15, 23.84) -- (134.15,105.93);

\path[draw=drawColor,line width= 0.4pt,line join=round,line cap=round] (136.01, 23.84) -- (136.01,121.63);

\path[draw=drawColor,line width= 0.4pt,line join=round,line cap=round] (137.86, 23.84) -- (137.86,127.01);

\path[draw=drawColor,line width= 0.4pt,line join=round,line cap=round] (139.72, 23.84) -- (139.72,132.40);

\path[draw=drawColor,line width= 0.4pt,line join=round,line cap=round] (141.58, 23.84) -- (141.58,117.59);

\path[draw=drawColor,line width= 0.4pt,line join=round,line cap=round] (143.44, 23.84) -- (143.44,109.52);

\path[draw=drawColor,line width= 0.4pt,line join=round,line cap=round] (145.30, 23.84) -- (145.30,101.00);

\path[draw=drawColor,line width= 0.4pt,line join=round,line cap=round] (147.15, 23.84) -- (147.15,110.87);

\path[draw=drawColor,line width= 0.4pt,line join=round,line cap=round] (149.01, 23.84) -- (149.01,111.76); 

\path[draw=drawColor,line width= 0.4pt,line join=round,line cap=round] (150.87, 23.84) -- (150.87,111.76);

\path[draw=drawColor,line width= 0.4pt,line join=round,line cap=round] (152.73, 23.84) -- (152.73,129.71);

\path[draw=drawColor,line width= 0.4pt,line join=round,line cap=round] (154.59, 23.84) -- (154.59,129.26);

\path[draw=drawColor,line width= 0.4pt,line join=round,line cap=round] (156.44, 23.84) -- (156.44,126.57);

\path[draw=drawColor,line width= 0.4pt,line join=round,line cap=round] (158.30, 23.84) -- (158.30,132.85);

\path[draw=drawColor,line width= 0.4pt,line join=round,line cap=round] (160.16, 23.84) -- (160.16,142.27);

\path[draw=drawColor,line width= 0.4pt,line join=round,line cap=round] (162.02, 23.84) -- (162.02,145.85);

\path[draw=drawColor,line width= 0.4pt,line join=round,line cap=round] (163.88, 23.84) -- (163.88,130.15);

\path[draw=drawColor,line width= 0.4pt,line join=round,line cap=round] (165.73, 23.84) -- (165.73,118.49);

\path[draw=drawColor,line width= 0.4pt,line join=round,line cap=round] (167.59, 23.84) -- (167.59,104.59);

\path[draw=drawColor,line width= 0.4pt,line join=round,line cap=round] (169.45, 23.84) -- (169.45,114.01);

\path[draw=drawColor,line width= 0.4pt,line join=round,line cap=round] (171.31, 23.84) -- (171.31,115.35); 

\path[draw=drawColor,line width= 0.4pt,line join=round,line cap=round] (173.17, 23.84) -- (173.17,108.17);

\path[draw=drawColor,line width= 0.4pt,line join=round,line cap=round] (175.02, 23.84) -- (175.02,129.26);

\path[draw=drawColor,line width= 0.4pt,line join=round,line cap=round] (176.88, 23.84) -- (176.88,125.67);

\path[draw=drawColor,line width= 0.4pt,line join=round,line cap=round] (178.74, 23.84) -- (178.74,128.81);

\path[draw=drawColor,line width= 0.4pt,line join=round,line cap=round] (180.60, 23.84) -- (180.60,142.27);

\path[draw=drawColor,line width= 0.4pt,line join=round,line cap=round] (182.46, 23.84) -- (182.46,159.31);

\path[draw=drawColor,line width= 0.4pt,line join=round,line cap=round] (184.31, 23.84) -- (184.31,155.27);

\path[draw=drawColor,line width= 0.4pt,line join=round,line cap=round] (186.17, 23.84) -- (186.17,140.02);

\path[draw=drawColor,line width= 0.4pt,line join=round,line cap=round] (188.03, 23.84) -- (188.03,126.57);

\path[draw=drawColor,line width= 0.4pt,line join=round,line cap=round] (189.89, 23.84) -- (189.89,114.90);

\path[draw=drawColor,line width= 0.4pt,line join=round,line cap=round] (191.75, 23.84) -- (191.75,126.57);

\path[draw=drawColor,line width= 0.4pt,line join=round,line cap=round] (193.60, 23.84) -- (193.60,132.40); 

\path[draw=drawColor,line width= 0.4pt,line join=round,line cap=round] (195.46, 23.84) -- (195.46,128.36);

\path[draw=drawColor,line width= 0.4pt,line join=round,line cap=round] (197.32, 23.84) -- (197.32,143.61);

\path[draw=drawColor,line width= 0.4pt,line join=round,line cap=round] (199.18, 23.84) -- (199.18,144.51);

\path[draw=drawColor,line width= 0.4pt,line join=round,line cap=round] (201.04, 23.84) -- (201.04,144.96);

\path[draw=drawColor,line width= 0.4pt,line join=round,line cap=round] (202.89, 23.84) -- (202.89,165.14);

\path[draw=drawColor,line width= 0.4pt,line join=round,line cap=round] (204.75, 23.84) -- (204.75,187.12);

\path[draw=drawColor,line width= 0.4pt,line join=round,line cap=round] (206.61, 23.84) -- (206.61,179.50);

\path[draw=drawColor,line width= 0.4pt,line join=round,line cap=round] (208.47, 23.84) -- (208.47,163.80);

\path[draw=drawColor,line width= 0.4pt,line join=round,line cap=round] (210.33, 23.84) -- (210.33,146.75);

\path[draw=drawColor,line width= 0.4pt,line join=round,line cap=round] (212.18, 23.84) -- (212.18,130.15);

\path[draw=drawColor,line width= 0.4pt,line join=round,line cap=round] (214.04, 23.84) -- (214.04,148.55);

\path[draw=drawColor,line width= 0.4pt,line join=round,line cap=round] (215.90, 23.84) -- (215.90,151.24); 

\path[draw=drawColor,line width= 0.4pt,line join=round,line cap=round] (217.76, 23.84) -- (217.76,148.10);

\path[draw=drawColor,line width= 0.4pt,line join=round,line cap=round] (219.62, 23.84) -- (219.62,166.04);

\path[draw=drawColor,line width= 0.4pt,line join=round,line cap=round] (221.47, 23.84) -- (221.47,164.25);

\path[draw=drawColor,line width= 0.4pt,line join=round,line cap=round] (223.33, 23.84) -- (223.33,166.49);

\path[draw=drawColor,line width= 0.4pt,line join=round,line cap=round] (225.19, 23.84) -- (225.19,191.61);

\path[draw=drawColor,line width= 0.4pt,line join=round,line cap=round] (227.05, 23.84) -- (227.05,209.10);

\path[draw=drawColor,line width= 0.4pt,line join=round,line cap=round] (228.91, 23.84) -- (228.91,205.52);

\path[draw=drawColor,line width= 0.4pt,line join=round,line cap=round] (230.76, 23.84) -- (230.76,183.09);

\path[draw=drawColor,line width= 0.4pt,line join=round,line cap=round] (232.62, 23.84) -- (232.62,161.11);

\path[draw=drawColor,line width= 0.4pt,line join=round,line cap=round] (234.48, 23.84) -- (234.48,145.41);

\path[draw=drawColor,line width= 0.4pt,line join=round,line cap=round] (236.34, 23.84) -- (236.34,161.11);

\path[draw=drawColor,line width= 0.4pt,line join=round,line cap=round] (238.20, 23.84) -- (238.20,165.14); 

\path[draw=drawColor,line width= 0.4pt,line join=round,line cap=round] (240.05, 23.84) -- (240.05,158.86);

\path[draw=drawColor,line width= 0.4pt,line join=round,line cap=round] (241.91, 23.84) -- (241.91,183.53);

\path[draw=drawColor,line width= 0.4pt,line join=round,line cap=round] (243.77, 23.84) -- (243.77,179.95);

\path[draw=drawColor,line width= 0.4pt,line join=round,line cap=round] (245.63, 23.84) -- (245.63,183.09);

\path[draw=drawColor,line width= 0.4pt,line join=round,line cap=round] (247.49, 23.84) -- (247.49,213.14);

\path[draw=drawColor,line width= 0.4pt,line join=round,line cap=round] (249.34, 23.84) -- (249.34,232.43);

\path[draw=drawColor,line width= 0.4pt,line join=round,line cap=round] (251.20, 23.84) -- (251.20,233.33);

\path[draw=drawColor,line width= 0.4pt,line join=round,line cap=round] (253.06, 23.84) -- (253.06,205.07);

\path[draw=drawColor,line width= 0.4pt,line join=round,line cap=round] (254.92, 23.84) -- (254.92,179.50);

\path[draw=drawColor,line width= 0.4pt,line join=round,line cap=round] (256.78, 23.84) -- (256.78,160.66);

\path[draw=drawColor,line width= 0.4pt,line join=round,line cap=round] (258.63, 23.84) -- (258.63,174.56);

\path[draw=drawColor,line width= 0.4pt,line join=round,line cap=round] (260.49, 23.84) -- (260.49,176.36); 

\path[draw=drawColor,line width= 0.4pt,line join=round,line cap=round] (262.35, 23.84) -- (262.35,166.49);

\path[draw=drawColor,line width= 0.4pt,line join=round,line cap=round] (264.21, 23.84) -- (264.21,186.23);

\path[draw=drawColor,line width= 0.4pt,line join=round,line cap=round] (266.07, 23.84) -- (266.07,179.95);

\path[draw=drawColor,line width= 0.4pt,line join=round,line cap=round] (267.92, 23.84) -- (267.92,186.67);

\path[draw=drawColor,line width= 0.4pt,line join=round,line cap=round] (269.78, 23.84) -- (269.78,218.97);

\path[draw=drawColor,line width= 0.4pt,line join=round,line cap=round] (271.64, 23.84) -- (271.64,244.09);

\path[draw=drawColor,line width= 0.4pt,line join=round,line cap=round] (273.50, 23.84) -- (273.50,250.37);

\path[draw=drawColor,line width= 0.4pt,line join=round,line cap=round] (275.36, 23.84) -- (275.36,205.07);

\path[draw=drawColor,line width= 0.4pt,line join=round,line cap=round] (277.21, 23.84) -- (277.21,184.88);

\path[draw=drawColor,line width= 0.4pt,line join=round,line cap=round] (279.07, 23.84) -- (279.07,162.90);

\path[draw=drawColor,line width= 0.4pt,line join=round,line cap=round] (280.93, 23.84) -- (280.93,175.01);

\path[draw=drawColor,line width= 0.4pt,line join=round,line cap=round] (282.79, 23.84) -- (282.79,185.33); 

\path[draw=drawColor,line width= 0.4pt,line join=round,line cap=round] (284.65, 23.84) -- (284.65,177.25);

\path[draw=drawColor,line width= 0.4pt,line join=round,line cap=round] (286.50, 23.84) -- (286.50,205.96);

\path[draw=drawColor,line width= 0.4pt,line join=round,line cap=round] (288.36, 23.84) -- (288.36,201.48);

\path[draw=drawColor,line width= 0.4pt,line join=round,line cap=round] (290.22, 23.84) -- (290.22,212.24);

\path[draw=drawColor,line width= 0.4pt,line join=round,line cap=round] (292.08, 23.84) -- (292.08,235.57);

\path[draw=drawColor,line width= 0.4pt,line join=round,line cap=round] (293.94, 23.84) -- (293.94,269.66);

\path[draw=drawColor,line width= 0.4pt,line join=round,line cap=round] (295.79, 23.84) -- (295.79,274.60);

\path[draw=drawColor,line width= 0.4pt,line join=round,line cap=round] (297.65, 23.84) -- (297.65,231.53);

\path[draw=drawColor,line width= 0.4pt,line join=round,line cap=round] (299.51, 23.84) -- (299.51,206.41);

\path[draw=drawColor,line width= 0.4pt,line join=round,line cap=round] (301.37, 23.84) -- (301.37,186.23);

\path[draw=drawColor,line width= 0.4pt,line join=round,line cap=round] (303.23, 23.84) -- (303.23,205.52);

\path[draw=drawColor,line width= 0.4pt,line join=round,line cap=round] (305.08, 23.84) -- (305.08,210.90); 

\path[draw=drawColor,line width= 0.4pt,line join=round,line cap=round] (306.94, 23.84) -- (306.94,199.24);

\path[draw=drawColor,line width= 0.4pt,line join=round,line cap=round] (308.80, 23.84) -- (308.80,211.80);

\path[draw=drawColor,line width= 0.4pt,line join=round,line cap=round] (310.66, 23.84) -- (310.66,230.64);

\path[draw=drawColor,line width= 0.4pt,line join=round,line cap=round] (312.52, 23.84) -- (312.52,235.57);

\path[draw=drawColor,line width= 0.4pt,line join=round,line cap=round] (314.37, 23.84) -- (314.37,263.83);

\path[draw=drawColor,line width= 0.4pt,line join=round,line cap=round] (316.23, 23.84) -- (316.23,302.86);

\path[draw=drawColor,line width= 0.4pt,line join=round,line cap=round] (318.09, 23.84) -- (318.09,295.68);

\path[draw=drawColor,line width= 0.4pt,line join=round,line cap=round] (319.95, 23.84) -- (319.95,251.72);

\path[draw=drawColor,line width= 0.4pt,line join=round,line cap=round] (321.81, 23.84) -- (321.81,230.64);

\path[draw=drawColor,line width= 0.4pt,line join=round,line cap=round] (323.66, 23.84) -- (323.66,198.79);

\path[draw=drawColor,line width= 0.4pt,line join=round,line cap=round] (325.52, 23.84) -- (325.52,217.63);
\end{scope}
\end{tikzpicture}
\end{center}
\caption{A time series: A discrete time stochastic process with 144 values corresponding to the number of 
   airlines passengers (in thousands) between 1949 and 1960.}\label{AirPassengers}
\end{figure}
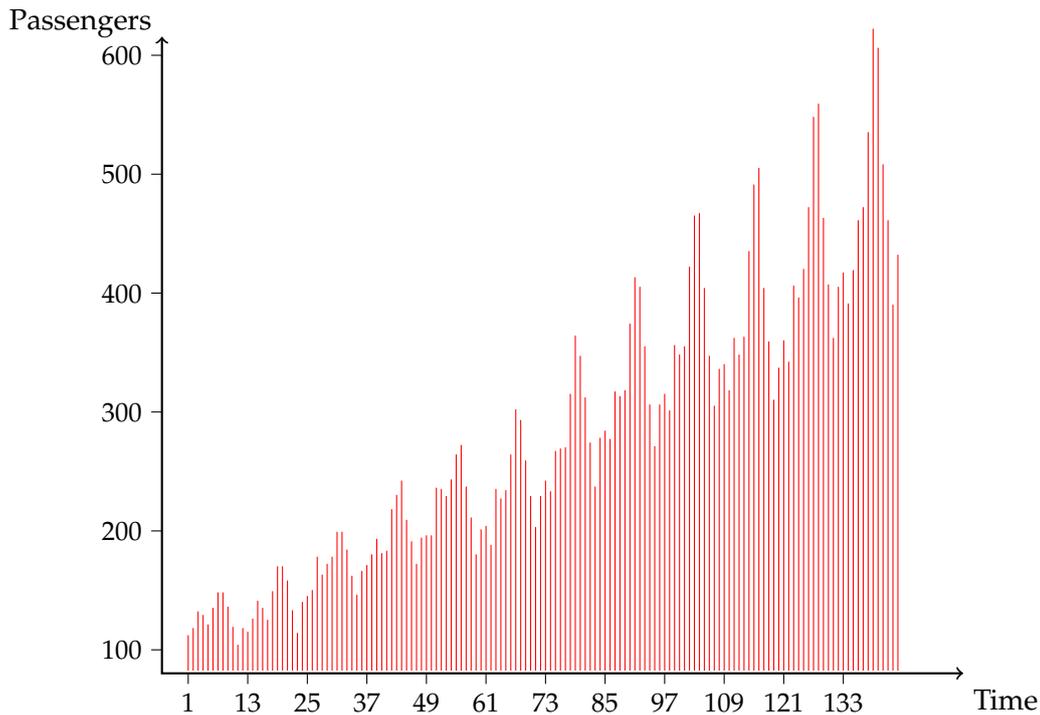

  Another important class of stochastic process in discrete time is the
  \Conc{time series}, that models a   
  different sort of dependency between variables. Figure \ref{AirPassengers} shows the monthly evolution
  of the number of passengers of international airlines between January 1949 and December 1960. One observes
  a \emph{trend} (increasing), a \emph{seasonality} (peaks at the central months of the year) and a residual 
  \emph{noise}
  (the purely random component of the process). Usually, one tries to fit a suitable model of dependence between
  the variables, so that the original process is expressed as the sum of these individual components.

\section{Example from biology: genes}
\begin{figure}[h]
  \begin{center}  
    \begin{overpic}
     [scale=2,tics=10]{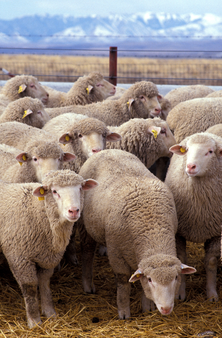}
    \end{overpic}
  \end{center}
  \caption{Some sheep}\label{somesheep}
\end{figure}

\begin{figure}[t]
  \begin{center}  
    \begin{overpic}
     [scale=0.5,tics=10]{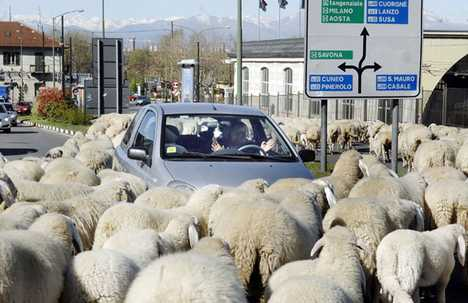}
    \end{overpic}
  \end{center}
  \caption{Some more sheep!}\label{somemoresheep}
\end{figure}

\subsection{Genotype and gene frequencies} 
    \Conc{Alleles} are several forms that a gene may have in a particular place 
    (\Conc{locus)} of a chromosome.

  For example, sheep haemoglobin presents two forms, produced by two alleles, $A$ and $B$, 
  of a certain locus.  
    Each individual possesses chromosomes in pairs, one coming from each parent.  
This implies that there are three possible \Conc{genotypes}: $AA$, $AB$, $BB$. 

    Typically, one allele is \Conc{dominant}, while the other is \Conc{recessive}. The recessive 
    allele shows up externally in the \Conc{phenotype} only if the dominant is not present.

    Assume we extract blood from a population of $N$ sheep, and the genotypes appear in proportions  
    $P_{AA}$, $P_{AB}$ and $P_{BB}$ (called \Conc{genotypic frequencies}).
    The \Conc{gene frequencies} are the proportions of the two alleles:     
    \begin{align}\label{eq:genefreq}
    \begin{split}
    P_A &:= P_{AA} + \frac{1}{2} P_{AB}
    \\
    P_B &:= P_{BB} + \frac{1}{2} P_{AB}
    \end{split}
    \end{align}

\subsection{Hardy-Weinberg principle}\label{sec:HWprinciple}
    Assume that:
    \begin{itemize}
    \item
      The proportions are the same for males and females. 
    \item  
      The genotype does not influence mating preferences.
    \item  
      Each allele of a parent is chosen with equal probability 1/2.
    \end{itemize}
    
  Then, the probabilities of each mating are, approximately (assuming a large population):
    \begin{align*}
    P(\text{$AA$ with $AA$}) &= P_{AA}^2
    \\
    P(\text{$AB$ with $AB$}) &= P_{AB}^2
    \\
    P(\text{$BB$ with $BB$}) &= P_{BB}^2
    \\
    P(\text{$AA$ with $AB$}) &= 2P_{AA}P_{AB}
    \\
    P(\text{$AA$ with $BB$}) &= 2P_{AA}P_{BB}
    \\
    P(\text{$AB$ with $BB$}) &= 2P_{AB}P_{BB}
    \end{align*}

  We can deduce easily the law of the genotypes for the next generation:
    \begin{align*}
    Q_{AA} &= P_{AA}^2 + \frac{1}{2}2P_{AA}P_{AB}+\frac{1}{4}P_{AB}^2 = P_A^2
    \\
    Q_{BB} &= P_{BB}^2 + \frac{1}{2}2P_{BB}P_{AB}+\frac{1}{4}P_{AB}^2 = P_B^2
    \\
    Q_{AB} &= 2P_AP_B
    \end{align*}
    
  Computing the gene frequencies $Q_A$ and $Q_B$ with (\ref{eq:genefreq}) we find again $P_A$ and $P_B$, 
  so the genotype frequencies  
  must be constant from the first generation onwards. This is the \Conc{Hardy-Weinberg} 
  principle (1908).
  
  As an application of this principle, suppose $B$ is recessive and we observe a 4\% proportion of individuals 
  showing the corresponding phenotype. Then we can deduce the genotype proportions
  of the whole population:
\begin{equation*}
  4\% = P_{BB}=P_B^2,\ \Rightarrow P_B=20\%,\ P_A=80\%, P_{AA}=64\%, P_{AB}=32\%
  \ .
\end{equation*}  
  If the population were small, then randomness in the mating may lead to \Conc{genetic drift}, and eventually
  one of the alleles will disappear from the population. The other gets \emph{fixed}, and time to fixation is one of
  the typical things of interest. This purely random fact explains the lost of genetic diversity in closed 
  small populations.

\subsection{Wright-Fisher model (1931)}\label{sec:WrightFisher}  

    If the mating is completely random, 
    it does not matter how the alleles are distributed among
    the $N$ individuals. We can simply consider the population of $2N$ alleles. 
    
    Assume that at generation 0 there are $X_0$ alleles of type $A$, 
    with $0<X_0<2N$. We pick alleles from this population independently from each other $2N$ times to form 
    the $N$ individuals of generation 1. 
    
  The law of the number of alleles of type $A$ must be $\Binom(2N,p)$, with $p=X_0/2N$. Thus
  \begin{align*}
    P\{X_1=k\}={2N\choose k} (X_0/2N)^k (1- X_0/2N)^{2N-k}
    \ ,\quad k=0,\dots, 2N
    \ .
  \end{align*}  
  In general, the number of alleles $A$ in generation $n+1$ knowing that there are $j$ in 
  generation $n$ is
  \begin{align*}
    P\{\condprob{X_{n+1}=k}{X_n=j}\}={2N\choose k} (j/2N)^k (1- j/2N)^{2N-k}
    \ .
  \end{align*}  
    This defines a Markov chain in discrete time. Its expectation is constant, $\E[X_n]=\E[X_0]$,
    and the expectation of the random variable $X_n$, knowing that a past variable $X_m$ ($m<n$)
    has taken value $k$, is equal to $k$:
  \begin{align}\label{eq:dtmartingale}
  \E\big[\condprob{X_{n}}{X_m=k}\big]=k
  \ ,\quad k=0,\dots,2N
  \ .
  \end{align}  
 
  However, as we saw in Section \ref{sec:HWprinciple}, the process will eventually reach states 0 or 2N, 
  and it will remain there forever. They are absorbing states (see Section \ref{longrunMC}).
     
\subsection{Conditional expectation}\label{sec:condexp}   
  The expression on the right-hand side of (\ref{eq:dtmartingale}) is called the \Conc{conditional expectation} of      
    $X_n$ given that $X_m=k$. It is exactly the expectation of $X_n$ computed from the conditional probability
    to the event $\{X_m=k\}$. One may write (\ref{eq:dtmartingale}) as
  \begin{align*}\label{eq:dtmartingale2}
  \E\big[\condprob{X_{n}}{X_m}\big]=X_m
  \end{align*}  
  and the left-hand side is now a random variable instead of a single number, called the \Conc{conditional
  expectation} of $X_n$ given $X_m$. For each $\omega\in\Omega$, the random variable $\E\big[\condprob{X}{Y}\big](\omega)$
  is equal to the number $\E\big[\condprob{X}{Y=y}\big]$, if $Y(\omega)=y$. 
  
  In case the conditioning random variable $Y$ has a continuous law, the definition above does not work, since 
  $\{Y(\omega)=y\}$ is an event of probability zero. The intuitive meaning is however the same. Mathematically, the trick 
  is not to consider the $y$ individually, but 
  collectively: The conditional expectation
  $\E\big[\condprob{X}{Y}\big]$ is given a sense as the (unique) random variable
  that can be factorized as a composition $(\varphi\circ Y)(\omega)$, with $\varphi\colon\Reals\rightarrow\Reals$, and 
  whose expectation, restricted to the events of $Y$, coincides with that of $X$:
 \begin{equation*}
  \E[(\varphi\circ Y)\cdot \ind_{\{Y\in B\}}] = \E[X\cdot \ind_{\{Y\in B\}} ]
  \ , 
\end{equation*}   
  where $\ind_{\{Y\in B\}}$ is equal to 1 if $Y(\omega)\in B$ and 0 otherwise.
\subsection{Continuous approximation of discrete laws} 


  Discrete laws involve only elementary discrete mathematics, but they are sometimes 
  cumbersome with computations. For instance, computing exactly the probability density 
  of a $\Binom(n,p)$ distribution when $n$ is large involve making the computer work
  with higher precision than usual. Although nowadays this is not a big deal (unless $n$
  is really very large), it is still useful, and conceptually important, to use continuous
  laws as a proxy to the real distribution.
  
\begin{figure}
  \begin{center}
\begin{tikzpicture}[x=1pt,y=1pt]
\definecolor{fillColor}{RGB}{255,255,255}
\path[use as bounding box,fill=fillColor,fill opacity=0.00] (50,50) rectangle (311.35,311.35);
\begin{scope}
\path[clip] ( 54.12, 67.32) rectangle (333.63,307.23);
\definecolor{drawColor}{RGB}{0,0,0}

\path[draw=drawColor,line width= 0.4pt,line join=round,line cap=round] (  4.45, 76.21) -- (  4.45, 76.21);

\path[draw=drawColor,line width= 0.4pt,line join=round,line cap=round] ( 12.87, 76.21) -- ( 12.87, 76.21);

\path[draw=drawColor,line width= 0.4pt,line join=round,line cap=round] ( 21.29, 76.21) -- ( 21.29, 76.21);

\path[draw=drawColor,line width= 0.4pt,line join=round,line cap=round] ( 29.71, 76.21) -- ( 29.71, 76.21);

\path[draw=drawColor,line width= 0.4pt,line join=round,line cap=round] ( 38.12, 76.21) -- ( 38.12, 76.21);

\path[draw=drawColor,line width= 0.4pt,line join=round,line cap=round] ( 46.54, 76.21) -- ( 46.54, 76.21);

\path[draw=drawColor,line width= 0.4pt,line join=round,line cap=round] ( 54.96, 76.21) -- ( 54.96, 76.21);

\path[draw=drawColor,line width= 0.4pt,line join=round,line cap=round] ( 63.38, 76.21) -- ( 63.38, 76.21);

\path[draw=drawColor,line width= 0.4pt,line join=round,line cap=round] ( 71.80, 76.21) -- ( 71.80, 76.21);

\path[draw=drawColor,line width= 0.4pt,line join=round,line cap=round] ( 80.22, 76.21) -- ( 80.22, 76.23);

\path[draw=drawColor,line width= 0.4pt,line join=round,line cap=round] ( 88.64, 76.21) -- ( 88.64, 76.30);

\path[draw=drawColor,line width= 0.4pt,line join=round,line cap=round] ( 97.06, 76.21) -- ( 97.06, 76.50);

\path[draw=drawColor,line width= 0.4pt,line join=round,line cap=round] (105.48, 76.21) -- (105.48, 77.02);

\path[draw=drawColor,line width= 0.4pt,line join=round,line cap=round] (113.90, 76.21) -- (113.90, 78.28);

\path[draw=drawColor,line width= 0.4pt,line join=round,line cap=round] (122.31, 76.21) -- (122.31, 81.03);

\path[draw=drawColor,line width= 0.4pt,line join=round,line cap=round] (130.73, 76.21) -- (130.73, 86.53);

\path[draw=drawColor,line width= 0.4pt,line join=round,line cap=round] (139.15, 76.21) -- (139.15, 96.45);

\path[draw=drawColor,line width= 0.4pt,line join=round,line cap=round] (147.57, 76.21) -- (147.57,112.54);

\path[draw=drawColor,line width= 0.4pt,line join=round,line cap=round] (155.99, 76.21) -- (155.99,135.87);

\path[draw=drawColor,line width= 0.4pt,line join=round,line cap=round] (164.41, 76.21) -- (164.41,165.70);

\path[draw=drawColor,line width= 0.4pt,line join=round,line cap=round] (172.83, 76.21) -- (172.83,198.61);

\path[draw=drawColor,line width= 0.4pt,line join=round,line cap=round] (181.25, 76.21) -- (181.25,228.53);

\path[draw=drawColor,line width= 0.4pt,line join=round,line cap=round] (189.67, 76.21) -- (189.67,248.19);

\path[draw=drawColor,line width= 0.4pt,line join=round,line cap=round] (198.08, 76.21) -- (198.08,251.77);

\path[draw=drawColor,line width= 0.4pt,line join=round,line cap=round] (206.50, 76.21) -- (206.50,237.59);

\path[draw=drawColor,line width= 0.4pt,line join=round,line cap=round] (214.92, 76.21) -- (214.92,209.11);

\path[draw=drawColor,line width= 0.4pt,line join=round,line cap=round] (223.34, 76.21) -- (223.34,173.67);

\path[draw=drawColor,line width= 0.4pt,line join=round,line cap=round] (231.76, 76.21) -- (231.76,139.37);

\path[draw=drawColor,line width= 0.4pt,line join=round,line cap=round] (240.18, 76.21) -- (240.18,112.06);

\path[draw=drawColor,line width= 0.4pt,line join=round,line cap=round] (248.60, 76.21) -- (248.60, 93.82);

\path[draw=drawColor,line width= 0.4pt,line join=round,line cap=round] (257.02, 76.21) -- (257.02, 83.58);

\path[draw=drawColor,line width= 0.4pt,line join=round,line cap=round] (265.44, 76.21) -- (265.44, 78.79);

\path[draw=drawColor,line width= 0.4pt,line join=round,line cap=round] (273.85, 76.21) -- (273.85, 76.94);

\path[draw=drawColor,line width= 0.4pt,line join=round,line cap=round] (282.27, 76.21) -- (282.27, 76.37);

\path[draw=drawColor,line width= 0.4pt,line join=round,line cap=round] (290.69, 76.21) -- (290.69, 76.23);

\path[draw=drawColor,line width= 0.4pt,line join=round,line cap=round] (299.11, 76.21) -- (299.11, 76.21);

\path[draw=drawColor,line width= 0.4pt,line join=round,line cap=round] (307.53, 76.21) -- (307.53, 76.21);
\end{scope}
\begin{scope}
\definecolor{drawColor}{RGB}{0,0,0}
\draw [<->, thick] ( 54.12,275.23) -- ( 54.12, 74.32) -- (333.63, 74.32);

\draw [-] ( 90.35, 70.32) node [below] {$-4$} -- ( 90.35, 74.32); 
\draw [-] ( 90.35, 70.32) node [below] {$-4$} -- ( 90.35, 74.32); 
\draw [-] (  142.11, 70.32) node [below] {$-2$} -- ( 142.11, 74.32); 
\draw [-] ( 90.35, 70.32) node [below] {$-4$} -- ( 90.35, 74.32); 
\draw [-] ( 193.87, 70.32) node [below] {$0$} -- ( 193.87, 74.32); 
\draw [-] ( 90.35, 70.32) node [below] {$-4$} -- ( 90.35, 74.32); 
\draw [-] ( 245.64, 70.32) node [below] {$2$} -- ( 245.64, 74.32); 
\draw [-] ( 90.35, 70.32) node [below] {$-4$} -- ( 90.35, 74.32); 
\draw [-] ( 297.40, 70.32) node [below] {$4$} -- ( 297.40, 74.32); 

\draw [-] ( 50.12, 76.21) node [left] {$0$} -- ( 54.12, 76.21);

\end{scope}
\begin{scope}
\path[clip] ( 54.12, 67.32) rectangle (333.63,307.23);
\definecolor{drawColor}{RGB}{255,0,0}
\path[draw=drawColor,line width= 0.4pt,line join=round,line cap=round] ( 64.47, 76.21) --
( 64.73, 76.21) --
( 64.99, 76.21) --
( 65.25, 76.21) --
( 65.51, 76.21) --
( 65.77, 76.21) --
( 66.03, 76.21) --
( 66.28, 76.21) --
( 66.54, 76.21) --
( 66.80, 76.21) --
( 67.06, 76.21) --
( 67.32, 76.21) --
( 67.58, 76.21) --
( 67.84, 76.21) --
( 68.10, 76.21) --
( 68.35, 76.21) --
( 68.61, 76.21) --
( 68.87, 76.21) --
( 69.13, 76.21) --
( 69.39, 76.21) --
( 69.65, 76.21) --
( 69.91, 76.21) --
( 70.17, 76.21) --
( 70.42, 76.21) --
( 70.68, 76.21) --
( 70.94, 76.21) --
( 71.20, 76.21) --
( 71.46, 76.21) --
( 71.72, 76.21) --
( 71.98, 76.21) --
( 72.24, 76.21) --
( 72.50, 76.21) --
( 72.75, 76.21) --
( 73.01, 76.21) --
( 73.27, 76.21) --
( 73.53, 76.21) --
( 73.79, 76.21) --
( 74.05, 76.21) --
( 74.31, 76.21) --
( 74.57, 76.21) --
( 74.82, 76.21) --
( 75.08, 76.21) --
( 75.34, 76.21) --
( 75.60, 76.21) --
( 75.86, 76.21) --
( 76.12, 76.21) --
( 76.38, 76.21) --
( 76.64, 76.21) --
( 76.89, 76.21) --
( 77.15, 76.21) --
( 77.41, 76.21) --
( 77.67, 76.21) --
( 77.93, 76.21) --
( 78.19, 76.21) --
( 78.45, 76.21) --
( 78.71, 76.21) --
( 78.97, 76.21) --
( 79.22, 76.22) --
( 79.48, 76.22) --
( 79.74, 76.22) --
( 80.00, 76.22) --
( 80.26, 76.22) --
( 80.52, 76.22) --
( 80.78, 76.22) --
( 81.04, 76.22) --
( 81.29, 76.22) --
( 81.55, 76.22) --
( 81.81, 76.22) --
( 82.07, 76.22) --
( 82.33, 76.22) --
( 82.59, 76.22) --
( 82.85, 76.22) --
( 83.11, 76.22) --
( 83.37, 76.23) --
( 83.62, 76.23) --
( 83.88, 76.23) --
( 84.14, 76.23) --
( 84.40, 76.23) --
( 84.66, 76.23) --
( 84.92, 76.23) --
( 85.18, 76.23) --
( 85.44, 76.23) --
( 85.69, 76.23) --
( 85.95, 76.24) --
( 86.21, 76.24) --
( 86.47, 76.24) --
( 86.73, 76.24) --
( 86.99, 76.24) --
( 87.25, 76.24) --
( 87.51, 76.24) --
( 87.76, 76.25) --
( 88.02, 76.25) --
( 88.28, 76.25) --
( 88.54, 76.25) --
( 88.80, 76.25) --
( 89.06, 76.25) --
( 89.32, 76.26) --
( 89.58, 76.26) --
( 89.84, 76.26) --
( 90.09, 76.26) --
( 90.35, 76.27) --
( 90.61, 76.27) --
( 90.87, 76.27) --
( 91.13, 76.27) --
( 91.39, 76.28) --
( 91.65, 76.28) --
( 91.91, 76.28) --
( 92.16, 76.28) --
( 92.42, 76.29) --
( 92.68, 76.29) --
( 92.94, 76.29) --
( 93.20, 76.30) --
( 93.46, 76.30) --
( 93.72, 76.30) --
( 93.98, 76.31) --
( 94.23, 76.31) --
( 94.49, 76.32) --
( 94.75, 76.32) --
( 95.01, 76.33) --
( 95.27, 76.33) --
( 95.53, 76.34) --
( 95.79, 76.34) --
( 96.05, 76.35) --
( 96.31, 76.35) --
( 96.56, 76.36) --
( 96.82, 76.36) --
( 97.08, 76.37) --
( 97.34, 76.37) --
( 97.60, 76.38) --
( 97.86, 76.39) --
( 98.12, 76.39) --
( 98.38, 76.40) --
( 98.63, 76.41) --
( 98.89, 76.42) --
( 99.15, 76.42) --
( 99.41, 76.43) --
( 99.67, 76.44) --
( 99.93, 76.45) --
(100.19, 76.46) --
(100.45, 76.47) --
(100.71, 76.48) --
(100.96, 76.49) --
(101.22, 76.50) --
(101.48, 76.51) --
(101.74, 76.52) --
(102.00, 76.53) --
(102.26, 76.54) --
(102.52, 76.55) --
(102.78, 76.57) --
(103.03, 76.58) --
(103.29, 76.59) --
(103.55, 76.61) --
(103.81, 76.62) --
(104.07, 76.64) --
(104.33, 76.65) --
(104.59, 76.67) --
(104.85, 76.68) --
(105.10, 76.70) --
(105.36, 76.72) --
(105.62, 76.73) --
(105.88, 76.75) --
(106.14, 76.77) --
(106.40, 76.79) --
(106.66, 76.81) --
(106.92, 76.83) --
(107.18, 76.85) --
(107.43, 76.88) --
(107.69, 76.90) --
(107.95, 76.92) --
(108.21, 76.95) --
(108.47, 76.97) --
(108.73, 77.00) --
(108.99, 77.02) --
(109.25, 77.05) --
(109.50, 77.08) --
(109.76, 77.11) --
(110.02, 77.14) --
(110.28, 77.17) --
(110.54, 77.20) --
(110.80, 77.23) --
(111.06, 77.26) --
(111.32, 77.30) --
(111.57, 77.33) --
(111.83, 77.37) --
(112.09, 77.41) --
(112.35, 77.45) --
(112.61, 77.49) --
(112.87, 77.53) --
(113.13, 77.57) --
(113.39, 77.61) --
(113.65, 77.66) --
(113.90, 77.70) --
(114.16, 77.75) --
(114.42, 77.80) --
(114.68, 77.85) --
(114.94, 77.90) --
(115.20, 77.95) --
(115.46, 78.00) --
(115.72, 78.06) --
(115.97, 78.12) --
(116.23, 78.17) --
(116.49, 78.23) --
(116.75, 78.30) --
(117.01, 78.36) --
(117.27, 78.42) --
(117.53, 78.49) --
(117.79, 78.56) --
(118.04, 78.63) --
(118.30, 78.70) --
(118.56, 78.77) --
(118.82, 78.85) --
(119.08, 78.93) --
(119.34, 79.01) --
(119.60, 79.09) --
(119.86, 79.17) --
(120.12, 79.26) --
(120.37, 79.35) --
(120.63, 79.44) --
(120.89, 79.53) --
(121.15, 79.62) --
(121.41, 79.72) --
(121.67, 79.82) --
(121.93, 79.92) --
(122.19, 80.03) --
(122.44, 80.14) --
(122.70, 80.25) --
(122.96, 80.36) --
(123.22, 80.47) --
(123.48, 80.59) --
(123.74, 80.71) --
(124.00, 80.84) --
(124.26, 80.96) --
(124.52, 81.09) --
(124.77, 81.22) --
(125.03, 81.36) --
(125.29, 81.50) --
(125.55, 81.64) --
(125.81, 81.78) --
(126.07, 81.93) --
(126.33, 82.08) --
(126.59, 82.24) --
(126.84, 82.40) --
(127.10, 82.56) --
(127.36, 82.73) --
(127.62, 82.90) --
(127.88, 83.07) --
(128.14, 83.25) --
(128.40, 83.43) --
(128.66, 83.61) --
(128.91, 83.80) --
(129.17, 83.99) --
(129.43, 84.19) --
(129.69, 84.39) --
(129.95, 84.60) --
(130.21, 84.81) --
(130.47, 85.02) --
(130.73, 85.24) --
(130.99, 85.46) --
(131.24, 85.69) --
(131.50, 85.92) --
(131.76, 86.15) --
(132.02, 86.40) --
(132.28, 86.64) --
(132.54, 86.89) --
(132.80, 87.15) --
(133.06, 87.41) --
(133.31, 87.68) --
(133.57, 87.95) --
(133.83, 88.22) --
(134.09, 88.50) --
(134.35, 88.79) --
(134.61, 89.08) --
(134.87, 89.38) --
(135.13, 89.68) --
(135.38, 89.99) --
(135.64, 90.31) --
(135.90, 90.63) --
(136.16, 90.95) --
(136.42, 91.28) --
(136.68, 91.62) --
(136.94, 91.97) --
(137.20, 92.32) --
(137.46, 92.67) --
(137.71, 93.03) --
(137.97, 93.40) --
(138.23, 93.78) --
(138.49, 94.16) --
(138.75, 94.55) --
(139.01, 94.94) --
(139.27, 95.34) --
(139.53, 95.75) --
(139.78, 96.16) --
(140.04, 96.58) --
(140.30, 97.01) --
(140.56, 97.44) --
(140.82, 97.88) --
(141.08, 98.33) --
(141.34, 98.79) --
(141.60, 99.25) --
(141.86, 99.72) --
(142.11,100.19) --
(142.37,100.68) --
(142.63,101.17) --
(142.89,101.66) --
(143.15,102.17) --
(143.41,102.68) --
(143.67,103.20) --
(143.93,103.73) --
(144.18,104.26) --
(144.44,104.81) --
(144.70,105.36) --
(144.96,105.91) --
(145.22,106.48) --
(145.48,107.05) --
(145.74,107.63) --
(146.00,108.22) --
(146.25,108.82) --
(146.51,109.42) --
(146.77,110.03) --
(147.03,110.65) --
(147.29,111.28) --
(147.55,111.92) --
(147.81,112.56) --
(148.07,113.21) --
(148.33,113.87) --
(148.58,114.54) --
(148.84,115.21) --
(149.10,115.89) --
(149.36,116.58) --
(149.62,117.28) --
(149.88,117.99) --
(150.14,118.70) --
(150.40,119.43) --
(150.65,120.16) --
(150.91,120.89) --
(151.17,121.64) --
(151.43,122.39) --
(151.69,123.15) --
(151.95,123.92) --
(152.21,124.70) --
(152.47,125.49) --
(152.72,126.28) --
(152.98,127.08) --
(153.24,127.89) --
(153.50,128.70) --
(153.76,129.52) --
(154.02,130.35) --
(154.28,131.19) --
(154.54,132.04) --
(154.80,132.89) --
(155.05,133.75) --
(155.31,134.61) --
(155.57,135.49) --
(155.83,136.37) --
(156.09,137.26) --
(156.35,138.15) --
(156.61,139.05) --
(156.87,139.96) --
(157.12,140.88) --
(157.38,141.80) --
(157.64,142.73) --
(157.90,143.66) --
(158.16,144.60) --
(158.42,145.55) --
(158.68,146.50) --
(158.94,147.46) --
(159.20,148.43) --
(159.45,149.40) --
(159.71,150.37) --
(159.97,151.35) --
(160.23,152.34) --
(160.49,153.33) --
(160.75,154.33) --
(161.01,155.33) --
(161.27,156.34) --
(161.52,157.35) --
(161.78,158.37) --
(162.04,159.39) --
(162.30,160.42) --
(162.56,161.44) --
(162.82,162.48) --
(163.08,163.52) --
(163.34,164.56) --
(163.59,165.60) --
(163.85,166.65) --
(164.11,167.70) --
(164.37,168.75) --
(164.63,169.81) --
(164.89,170.87) --
(165.15,171.93) --
(165.41,172.99) --
(165.67,174.06) --
(165.92,175.13) --
(166.18,176.19) --
(166.44,177.27) --
(166.70,178.34) --
(166.96,179.41) --
(167.22,180.48) --
(167.48,181.56) --
(167.74,182.63) --
(167.99,183.71) --
(168.25,184.78) --
(168.51,185.86) --
(168.77,186.93) --
(169.03,188.01) --
(169.29,189.08) --
(169.55,190.15) --
(169.81,191.22) --
(170.06,192.29) --
(170.32,193.36) --
(170.58,194.42) --
(170.84,195.48) --
(171.10,196.54) --
(171.36,197.60) --
(171.62,198.66) --
(171.88,199.71) --
(172.14,200.76) --
(172.39,201.80) --
(172.65,202.84) --
(172.91,203.88) --
(173.17,204.91) --
(173.43,205.94) --
(173.69,206.96) --
(173.95,207.98) --
(174.21,208.99) --
(174.46,209.99) --
(174.72,210.99) --
(174.98,211.99) --
(175.24,212.98) --
(175.50,213.96) --
(175.76,214.93) --
(176.02,215.90) --
(176.28,216.86) --
(176.54,217.81) --
(176.79,218.76) --
(177.05,219.69) --
(177.31,220.62) --
(177.57,221.54) --
(177.83,222.45) --
(178.09,223.36) --
(178.35,224.25) --
(178.61,225.13) --
(178.86,226.01) --
(179.12,226.87) --
(179.38,227.72) --
(179.64,228.57) --
(179.90,229.40) --
(180.16,230.22) --
(180.42,231.03) --
(180.68,231.83) --
(180.93,232.62) --
(181.19,233.40) --
(181.45,234.16) --
(181.71,234.91) --
(181.97,235.65) --
(182.23,236.38) --
(182.49,237.09) --
(182.75,237.80) --
(183.01,238.48) --
(183.26,239.16) --
(183.52,239.82) --
(183.78,240.47) --
(184.04,241.10) --
(184.30,241.72) --
(184.56,242.33) --
(184.82,242.92) --
(185.08,243.49) --
(185.33,244.05) --
(185.59,244.60) --
(185.85,245.13) --
(186.11,245.65) --
(186.37,246.15) --
(186.63,246.63) --
(186.89,247.10) --
(187.15,247.56) --
(187.40,247.99) --
(187.66,248.42) --
(187.92,248.82) --
(188.18,249.21) --
(188.44,249.58) --
(188.70,249.94) --
(188.96,250.28) --
(189.22,250.60) --
(189.48,250.90) --
(189.73,251.19) --
(189.99,251.46) --
(190.25,251.72) --
(190.51,251.96) --
(190.77,252.18) --
(191.03,252.38) --
(191.29,252.56) --
(191.55,252.73) --
(191.80,252.88) --
(192.06,253.01) --
(192.32,253.13) --
(192.58,253.23) --
(192.84,253.31) --
(193.10,253.37) --
(193.36,253.41) --
(193.62,253.44) --
(193.87,253.45) --
(194.13,253.44) --
(194.39,253.41) --
(194.65,253.37) --
(194.91,253.31) --
(195.17,253.23) --
(195.43,253.13) --
(195.69,253.01) --
(195.95,252.88) --
(196.20,252.73) --
(196.46,252.56) --
(196.72,252.38) --
(196.98,252.18) --
(197.24,251.96) --
(197.50,251.72) --
(197.76,251.46) --
(198.02,251.19) --
(198.27,250.90) --
(198.53,250.60) --
(198.79,250.28) --
(199.05,249.94) --
(199.31,249.58) --
(199.57,249.21) --
(199.83,248.82) --
(200.09,248.42) --
(200.35,247.99) --
(200.60,247.56) --
(200.86,247.10) --
(201.12,246.63) --
(201.38,246.15) --
(201.64,245.65) --
(201.90,245.13) --
(202.16,244.60) --
(202.42,244.05) --
(202.67,243.49) --
(202.93,242.92) --
(203.19,242.33) --
(203.45,241.72) --
(203.71,241.10) --
(203.97,240.47) --
(204.23,239.82) --
(204.49,239.16) --
(204.74,238.48) --
(205.00,237.80) --
(205.26,237.09) --
(205.52,236.38) --
(205.78,235.65) --
(206.04,234.91) --
(206.30,234.16) --
(206.56,233.40) --
(206.82,232.62) --
(207.07,231.83) --
(207.33,231.03) --
(207.59,230.22) --
(207.85,229.40) --
(208.11,228.57) --
(208.37,227.72) --
(208.63,226.87) --
(208.89,226.01) --
(209.14,225.13) --
(209.40,224.25) --
(209.66,223.36) --
(209.92,222.45) --
(210.18,221.54) --
(210.44,220.62) --
(210.70,219.69) --
(210.96,218.76) --
(211.21,217.81) --
(211.47,216.86) --
(211.73,215.90) --
(211.99,214.93) --
(212.25,213.96) --
(212.51,212.98) --
(212.77,211.99) --
(213.03,210.99) --
(213.29,209.99) --
(213.54,208.99) --
(213.80,207.98) --
(214.06,206.96) --
(214.32,205.94) --
(214.58,204.91) --
(214.84,203.88) --
(215.10,202.84) --
(215.36,201.80) --
(215.61,200.76) --
(215.87,199.71) --
(216.13,198.66) --
(216.39,197.60) --
(216.65,196.54) --
(216.91,195.48) --
(217.17,194.42) --
(217.43,193.36) --
(217.69,192.29) --
(217.94,191.22) --
(218.20,190.15) --
(218.46,189.08) --
(218.72,188.01) --
(218.98,186.93) --
(219.24,185.86) --
(219.50,184.78) --
(219.76,183.71) --
(220.01,182.63) --
(220.27,181.56) --
(220.53,180.48) --
(220.79,179.41) --
(221.05,178.34) --
(221.31,177.27) --
(221.57,176.19) --
(221.83,175.13) --
(222.08,174.06) --
(222.34,172.99) --
(222.60,171.93) --
(222.86,170.87) --
(223.12,169.81) --
(223.38,168.75) --
(223.64,167.70) --
(223.90,166.65) --
(224.16,165.60) --
(224.41,164.56) --
(224.67,163.52) --
(224.93,162.48) --
(225.19,161.44) --
(225.45,160.42) --
(225.71,159.39) --
(225.97,158.37) --
(226.23,157.35) --
(226.48,156.34) --
(226.74,155.33) --
(227.00,154.33) --
(227.26,153.33) --
(227.52,152.34) --
(227.78,151.35) --
(228.04,150.37) --
(228.30,149.40) --
(228.55,148.43) --
(228.81,147.46) --
(229.07,146.50) --
(229.33,145.55) --
(229.59,144.60) --
(229.85,143.66) --
(230.11,142.73) --
(230.37,141.80) --
(230.63,140.88) --
(230.88,139.96) --
(231.14,139.05) --
(231.40,138.15) --
(231.66,137.26) --
(231.92,136.37) --
(232.18,135.49) --
(232.44,134.61) --
(232.70,133.75) --
(232.95,132.89) --
(233.21,132.04) --
(233.47,131.19) --
(233.73,130.35) --
(233.99,129.52) --
(234.25,128.70) --
(234.51,127.89) --
(234.77,127.08) --
(235.03,126.28) --
(235.28,125.49) --
(235.54,124.70) --
(235.80,123.92) --
(236.06,123.15) --
(236.32,122.39) --
(236.58,121.64) --
(236.84,120.89) --
(237.10,120.16) --
(237.35,119.43) --
(237.61,118.70) --
(237.87,117.99) --
(238.13,117.28) --
(238.39,116.58) --
(238.65,115.89) --
(238.91,115.21) --
(239.17,114.54) --
(239.42,113.87) --
(239.68,113.21) --
(239.94,112.56) --
(240.20,111.92) --
(240.46,111.28) --
(240.72,110.65) --
(240.98,110.03) --
(241.24,109.42) --
(241.50,108.82) --
(241.75,108.22) --
(242.01,107.63) --
(242.27,107.05) --
(242.53,106.48) --
(242.79,105.91) --
(243.05,105.36) --
(243.31,104.81) --
(243.57,104.26) --
(243.82,103.73) --
(244.08,103.20) --
(244.34,102.68) --
(244.60,102.17) --
(244.86,101.66) --
(245.12,101.17) --
(245.38,100.68) --
(245.64,100.19) --
(245.89, 99.72) --
(246.15, 99.25) --
(246.41, 98.79) --
(246.67, 98.33) --
(246.93, 97.88) --
(247.19, 97.44) --
(247.45, 97.01) --
(247.71, 96.58) --
(247.97, 96.16) --
(248.22, 95.75) --
(248.48, 95.34) --
(248.74, 94.94) --
(249.00, 94.55) --
(249.26, 94.16) --
(249.52, 93.78) --
(249.78, 93.40) --
(250.04, 93.03) --
(250.29, 92.67) --
(250.55, 92.32) --
(250.81, 91.97) --
(251.07, 91.62) --
(251.33, 91.28) --
(251.59, 90.95) --
(251.85, 90.63) --
(252.11, 90.31) --
(252.37, 89.99) --
(252.62, 89.68) --
(252.88, 89.38) --
(253.14, 89.08) --
(253.40, 88.79) --
(253.66, 88.50) --
(253.92, 88.22) --
(254.18, 87.95) --
(254.44, 87.68) --
(254.69, 87.41) --
(254.95, 87.15) --
(255.21, 86.89) --
(255.47, 86.64) --
(255.73, 86.40) --
(255.99, 86.15) --
(256.25, 85.92) --
(256.51, 85.69) --
(256.76, 85.46) --
(257.02, 85.24) --
(257.28, 85.02) --
(257.54, 84.81) --
(257.80, 84.60) --
(258.06, 84.39) --
(258.32, 84.19) --
(258.58, 83.99) --
(258.84, 83.80) --
(259.09, 83.61) --
(259.35, 83.43) --
(259.61, 83.25) --
(259.87, 83.07) --
(260.13, 82.90) --
(260.39, 82.73) --
(260.65, 82.56) --
(260.91, 82.40) --
(261.16, 82.24) --
(261.42, 82.08) --
(261.68, 81.93) --
(261.94, 81.78) --
(262.20, 81.64) --
(262.46, 81.50) --
(262.72, 81.36) --
(262.98, 81.22) --
(263.23, 81.09) --
(263.49, 80.96) --
(263.75, 80.84) --
(264.01, 80.71) --
(264.27, 80.59) --
(264.53, 80.47) --
(264.79, 80.36) --
(265.05, 80.25) --
(265.31, 80.14) --
(265.56, 80.03) --
(265.82, 79.92) --
(266.08, 79.82) --
(266.34, 79.72) --
(266.60, 79.62) --
(266.86, 79.53) --
(267.12, 79.44) --
(267.38, 79.35) --
(267.63, 79.26) --
(267.89, 79.17) --
(268.15, 79.09) --
(268.41, 79.01) --
(268.67, 78.93) --
(268.93, 78.85) --
(269.19, 78.77) --
(269.45, 78.70) --
(269.71, 78.63) --
(269.96, 78.56) --
(270.22, 78.49) --
(270.48, 78.42) --
(270.74, 78.36) --
(271.00, 78.30) --
(271.26, 78.23) --
(271.52, 78.17) --
(271.78, 78.12) --
(272.03, 78.06) --
(272.29, 78.00) --
(272.55, 77.95) --
(272.81, 77.90) --
(273.07, 77.85) --
(273.33, 77.80) --
(273.59, 77.75) --
(273.85, 77.70) --
(274.10, 77.66) --
(274.36, 77.61) --
(274.62, 77.57) --
(274.88, 77.53) --
(275.14, 77.49) --
(275.40, 77.45) --
(275.66, 77.41) --
(275.92, 77.37) --
(276.18, 77.33) --
(276.43, 77.30) --
(276.69, 77.26) --
(276.95, 77.23) --
(277.21, 77.20) --
(277.47, 77.17) --
(277.73, 77.14) --
(277.99, 77.11) --
(278.25, 77.08) --
(278.50, 77.05) --
(278.76, 77.02) --
(279.02, 77.00) --
(279.28, 76.97) --
(279.54, 76.95) --
(279.80, 76.92) --
(280.06, 76.90) --
(280.32, 76.88) --
(280.57, 76.85) --
(280.83, 76.83) --
(281.09, 76.81) --
(281.35, 76.79) --
(281.61, 76.77) --
(281.87, 76.75) --
(282.13, 76.73) --
(282.39, 76.72) --
(282.65, 76.70) --
(282.90, 76.68) --
(283.16, 76.67) --
(283.42, 76.65) --
(283.68, 76.64) --
(283.94, 76.62) --
(284.20, 76.61) --
(284.46, 76.59) --
(284.72, 76.58) --
(284.97, 76.57) --
(285.23, 76.55) --
(285.49, 76.54) --
(285.75, 76.53) --
(286.01, 76.52) --
(286.27, 76.51) --
(286.53, 76.50) --
(286.79, 76.49) --
(287.04, 76.48) --
(287.30, 76.47) --
(287.56, 76.46) --
(287.82, 76.45) --
(288.08, 76.44) --
(288.34, 76.43) --
(288.60, 76.42) --
(288.86, 76.42) --
(289.12, 76.41) --
(289.37, 76.40) --
(289.63, 76.39) --
(289.89, 76.39) --
(290.15, 76.38) --
(290.41, 76.37) --
(290.67, 76.37) --
(290.93, 76.36) --
(291.19, 76.36) --
(291.44, 76.35) --
(291.70, 76.35) --
(291.96, 76.34) --
(292.22, 76.34) --
(292.48, 76.33) --
(292.74, 76.33) --
(293.00, 76.32) --
(293.26, 76.32) --
(293.52, 76.31) --
(293.77, 76.31) --
(294.03, 76.30) --
(294.29, 76.30) --
(294.55, 76.30) --
(294.81, 76.29) --
(295.07, 76.29) --
(295.33, 76.29) --
(295.59, 76.28) --
(295.84, 76.28) --
(296.10, 76.28) --
(296.36, 76.28) --
(296.62, 76.27) --
(296.88, 76.27) --
(297.14, 76.27) --
(297.40, 76.27) --
(297.66, 76.26) --
(297.91, 76.26) --
(298.17, 76.26) --
(298.43, 76.26) --
(298.69, 76.25) --
(298.95, 76.25) --
(299.21, 76.25) --
(299.47, 76.25) --
(299.73, 76.25) --
(299.99, 76.25) --
(300.24, 76.24) --
(300.50, 76.24) --
(300.76, 76.24) --
(301.02, 76.24) --
(301.28, 76.24) --
(301.54, 76.24) --
(301.80, 76.24) --
(302.06, 76.23) --
(302.31, 76.23) --
(302.57, 76.23) --
(302.83, 76.23) --
(303.09, 76.23) --
(303.35, 76.23) --
(303.61, 76.23) --
(303.87, 76.23) --
(304.13, 76.23) --
(304.38, 76.23) --
(304.64, 76.22) --
(304.90, 76.22) --
(305.16, 76.22) --
(305.42, 76.22) --
(305.68, 76.22) --
(305.94, 76.22) --
(306.20, 76.22) --
(306.46, 76.22) --
(306.71, 76.22) --
(306.97, 76.22) --
(307.23, 76.22) --
(307.49, 76.22) --
(307.75, 76.22) --
(308.01, 76.22) --
(308.27, 76.22) --
(308.53, 76.22) --
(308.78, 76.21) --
(309.04, 76.21) --
(309.30, 76.21) --
(309.56, 76.21) --
(309.82, 76.21) --
(310.08, 76.21) --
(310.34, 76.21) --
(310.60, 76.21) --
(310.86, 76.21) --
(311.11, 76.21) --
(311.37, 76.21) --
(311.63, 76.21) --
(311.89, 76.21) --
(312.15, 76.21) --
(312.41, 76.21) --
(312.67, 76.21) --
(312.93, 76.21) --
(313.18, 76.21) --
(313.44, 76.21) --
(313.70, 76.21) --
(313.96, 76.21) --
(314.22, 76.21) --
(314.48, 76.21) --
(314.74, 76.21) --
(315.00, 76.21) --
(315.25, 76.21) --
(315.51, 76.21) --
(315.77, 76.21) --
(316.03, 76.21) --
(316.29, 76.21) --
(316.55, 76.21) --
(316.81, 76.21) --
(317.07, 76.21) --
(317.33, 76.21) --
(317.58, 76.21) --
(317.84, 76.21) --
(318.10, 76.21) --
(318.36, 76.21) --
(318.62, 76.21) --
(318.88, 76.21) --
(319.14, 76.21) --
(319.40, 76.21) --
(319.65, 76.21) --
(319.91, 76.21) --
(320.17, 76.21) --
(320.43, 76.21) --
(320.69, 76.21) --
(320.95, 76.21) --
(321.21, 76.21) --
(321.47, 76.21) --
(321.72, 76.21) --
(321.98, 76.21) --
(322.24, 76.21) --
(322.50, 76.21) --
(322.76, 76.21) --
(323.02, 76.21) --
(323.28, 76.21);
\end{scope}
\end{tikzpicture}
\end{center}
\caption{Approximation by the Gaussian law: The probability density of $\Binom(n=45, p=0.7)$, after subtracting its expectation
  (\emph{centring}) and dividing by the square root of its variance (\emph{reducing}), depicted with black vertical lines.
  In red, the density function of $N(0,1)$.}\label{fig:BinomGauss}
\end{figure}
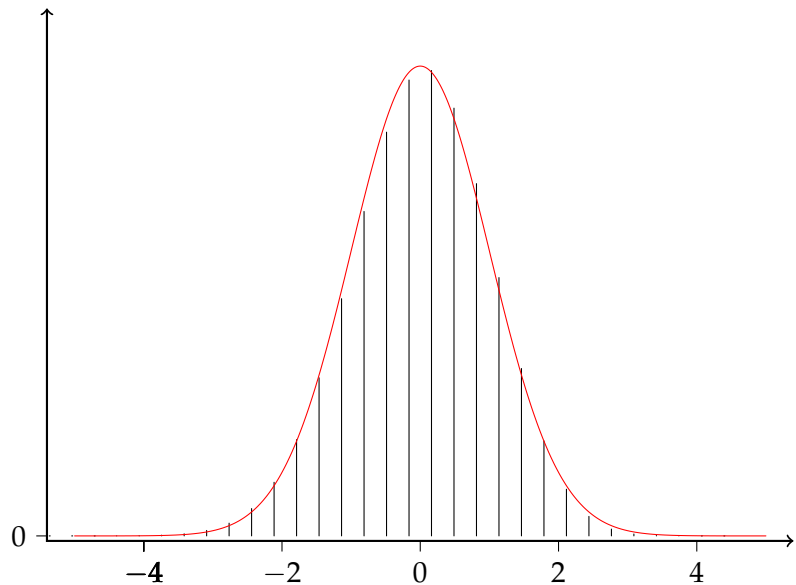

  Specifically, for the binomial law:
   If $X\sim \Binom(n,p)$, then 
  \begin{align*}
  \frac{X-np}{\sqrt{np(1-p)}}
  \sim 
  N(0,1)
  \ ,\quad \text{(approximately, for $n$ large)}
    \end{align*}  
  where   
    $N(0,1)$ denotes the so called \Conc{Normal} (or \Conc{Gaussian}) \Conc{law} with     
    expectation 0 and variance 1. Its density function is the \emph{Gaussian bell curve}
\begin{equation*}
    f(x)=\frac{1}{\sqrt{2\pi}}e^{-x^2/2}
    \ .
\end{equation*}   
  Figure \ref{fig:BinomGauss} shows graphically the approximation.

   The importance of the Gaussian law comes from the \Conc{Central Limit Theorem}, 
   which explains its ubiquity:
    If $X_n$ is a sequence of independent identically distributed random variables, 
    with finite variance, and $S_n:=\sum_{i=1}^n X_i$, then 
    \begin{align*}
    \frac{S_n - \E[S_n]}{\sqrt{\Var[S_n]}}
    \quad 
    \text{\emph{converges in law to}} 
    \quad
    N(0,1)
    \ .
    \end{align*}                   
    
  We immediately see that the binomial case above is a particular application of this theorem,
  taking $X_i\sim \Binom(1,p)$, which implies $S_n\sim\Binom(n,p)$.  
  Convergence in law is a non-elementary concept that has to do with  
  duality in functional spaces: Suppose that $\{Y_n\}$ is a sequence of
  random variables with respective distributions $P_n$, and $Y$ is a random variable with 
  distribution $P$. Then, we say that $\{Y_n\}$ \Conc{converge in law to $P$} if for every
  bounded continuous function $f\colon\Reals\rightarrow\Reals$, 
\begin{equation*}
  \lim_{n\to \infty}\E[f(Y_n)] = \E[f(Y)] 
  \ .
\end{equation*}  
  The seemingly natural ``setwise'' convergence $\lim_{n\to\infty} P_n(A) = P(A)$
  for all sets $A$
  is too strong, and will not work for the purpose of approximating by continuous distributions.
  
  One practical consequence of the Central Limit Theorem for modelling is that any phenomenon 
  whose result is the sum of many quantitatively small causes (like for instance the height
  or the weight of a person) will be well described by a Gaussian random variable. The fact that
  the Gaussian laws may take values in any interval of the real line is not an obstacle due to
  the rapid decrease of the bell curve: Outside a short interval, the probability is extremely small.
  
\subsection{Random walk and the Wiener process} 
    Let $\{X_n\}_n$ be a Markov chain taking values in $\Enters$ with  
    \begin{align*}
    X_0&=0
    \\
    P\{\condprob{X_{n+1}=i+1}{X_n=i}\}&=1/2
    \\
    P\{\condprob{X_{n+1}=i-1}{X_n=i}\}&=1/2
    \ .
    \end{align*}
   This process is called \Conc{random walk}. It is simply a ``walk'' on the integer lattice,
   where at each time step we go to the left or to right according to the toss of a fair coin.
   In other words, the increments $\varepsilon_n:=X_n-X_{n-1}$ are independent and take values 1 and -1 with probability 1/2.
    
    Define a sequence of continuous-time process $W^N_t$ by renormalisation of a random walk: 
  \begin{align*}
    W^N_t=\frac{1}{\sqrt{N}} X_{\lfloor Nt \rfloor} 
    =\frac{1}{\sqrt{N}} \sum_{k=1}^{\lfloor Nt \rfloor} \varepsilon_k 
    \ . 
  \end{align*}
  By the Central Limit Theorem, 
  \begin{align*}
  =\frac{1}{\sqrt{Nt}} \sum_{k=1}^{\lfloor Nt \rfloor} \varepsilon_k  
       \quad 
       \text{converges in law to} 
       \quad
       N(0,1)
       \ ,
  \end{align*}
  hence the sequence $\{W^N_t\}_N$ converges in law to a random variable $W_t\sim N(0,t)$, 
  the Gaussian law with variance $t$, for all $t>0$, whose density is 
  \begin{equation*}
  f(x)=\frac{1}{\sqrt{2\pi t}}e^{-x^2/2t}
  \ .
  \end{equation*}   
  Analogously, $W^N_t-W^N_s$ converges in law to $W_t-W_s \sim N(0,t-s)$. 
    The limiting process $W_t$ satisfies: 
    \begin{enumerate}
    \item
    The increments in non-overlapping intervals are independent.
    \item
    The expectation is constant equal to zero.
    \item
    The sample paths are continuous functions.
    \item
    The sample paths are non-differentiable at any point.
    \end{enumerate}
    
    $W$ is called the \Conc{Wiener process}. 
    In fact, a Wiener process is defined by its laws, but usually it is additionally
    asked to have continuous paths. This particular construction as the limit of 
    random walks 
    leads indeed to continuous paths.
    
    The Wiener process is also called \Conc{Brownian Motion} in the mathematical
    literature. However, the Brownian motion is a physical phenomenon, and the 
    Wiener process is just a mathematical model (and not the best one) to 
    that phenomenon.

\subsection{Diffusion approximation of Wright-Fisher model} 
   The Markov chain of the Wright-Fisher model is too
   complicated to work upon. Instead,
    define $Y^N_t = \frac{1}{2N} X^N_{\lfloor 2Nt \rfloor}$. Then,
    \begin{align*}
     &\E\big[\condprob{(Y^N_{t+h} - Y^N_t)^2}{Y^N_t=x}\big]  
     \\
     &=
      \Big(\frac{1}{2N}\Big)^2\E\big[\condprob{(X^N_{\lfloor 2N(t+h) \rfloor}-X^N_t)^2}{X^N_t=x}\big] 
     \\ 
     &= 
     \frac{1}{2N}x(1-x)\ ,\quad \text{(if $h\sim \frac{1}{2N}$)}.
    \end{align*}

    The limiting process $Y_t$ exists, satisfies
    \begin{align*}
     \E\big[\condprob{(Y_{t+h} - Y_t)^2}{Y_t=x}\big]  
     = 
     hx(1-x) +o(h)
    \end{align*}
    and it is called the \Conc{diffusion approximation} of the original Markov chain.
  
\subsection{Diffusions}\label{sec:diffusions}
    A \Conc{diffusion} $Y$ is a continuous-time Markov process, with continuous paths, and such
    that
    \begin{enumerate}
      \item 
        $\E[\condprob{Y_{t+h}-Y_t}{Y_t=x}] = b(t,x)h + o(h)$
      \item   
        $\E[\condprob{(Y_{t+h}-Y_t)^2}{Y_t=x}] = a(t,x)h + o(h)$
    \end{enumerate}
    for some functions $a$ and $b$. See Section \ref{sec:condexp}  for the interpretation of the conditional expectations
    when the conditioning variable is continuous.
    
    Under mild conditions, $Y_t$ has  
    a continuous law with density $f(t,x)$ satisfying
    the \Conc{Kolmogorov forward and backward equations}: 
    \begin{align}\label{eq:Kolmof}
        \displaystyle \frac{\partial}{\partial t}f(t,x) &= 
        \frac{1}{2}\frac{\partial^2}{\partial x^2} \big[a(t,x)f(t,x)]
        -\frac{\partial}{\partial x} \big[b(t,x)f(t,x)] 
        \\ \label{eq:Kolmob}
        \displaystyle -\frac{\partial}{\partial s}f(t,x) &= 
        \frac{1}{2}a(t,x)\frac{\partial^2}{\partial x^2} f(t,x)
        -b(t,x)\frac{\partial}{\partial x} f(t,x) 
    \end{align}

  The Wright--Fisher model can be expanded to take into account other effects in population dynamics,
  such as selection or mutation. This complications make even more useful the corresponding 
  diffusion approximations.
  
\section{Example from economy: stock markets}

\begin{figure}[h]
  \begin{center}  
    \begin{overpic}
     [scale=0.5,tics=10]{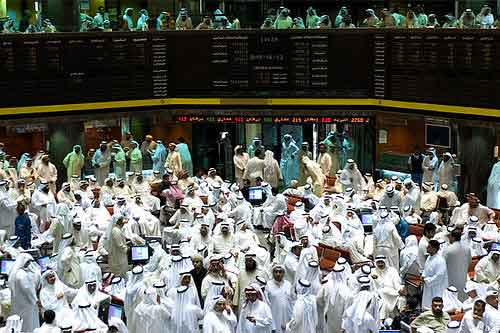}
    \end{overpic}
  \end{center}
  \caption{Kuwait stock market}\label{stockmarket}
\end{figure}

\subsection{A binomial economy} 
    Assume an economy with only two states: 
    \begin{itemize}
    \item
       Up (with probability $p$)
    \item   
       Down (with probability $1-p$)
    \end{itemize}
    
    Assume that there are two assets:
    \begin{itemize}
    \item
       A \Conc{risk-free bond} with interest rate $R$, and
    \item   
       A \Conc{share} with price $S(0)$ at time 0 and $S(1)$ at time 1, given
       by 
       \begin{align*}
       S(1)=
       \begin{cases}
         S(0)u\ ,\quad \text{if the economy is ``up''}
         \\
         S(0)d\ ,\quad \text{if the economy is ``down''}
       \end{cases}
       \end{align*}
    \end{itemize}
    
    A \Conc{trading strategy} for a portfolio is defined by 
    \begin{itemize}
    \item
      $B_0$ \euro\ allocated to the bond, and
    \item
      $\Delta_0$ quantity of shares of the stock
    \end{itemize}
    at time zero.
  The values of the portfolio at times 0 and 1 are
  \begin{align*}
    V(0)&=B_0+\Delta_0S(0)
    \\
    V(1)&=B_0(1+R)+\Delta_0S(1)
  \end{align*}

\subsection{Free lunch?} 
    As we will see, one can make money for free, unless $d< 1+R < u$.
    An \Conc{arbitrage opportunity} is the situation in which, without investing any 
    money at time zero, the probability
    to have a positive portfolio at time one is positive, 
    and the probability of a loss is zero.
    
   For a couple $(B_0,\Delta_0)$ 
   such that $V(0)=0$,
       \begin{align*}
       V(1)
       &= B_0(1+R)+\Delta_0S(1) 
       =
       \begin{cases}
         B_0(1+R)+\Delta_0uS(0)
         \\
         B_0(1+R)+\Delta_0dS(0)
       \end{cases}
       \\
       &=
       \begin{cases}
         \Delta_0S(0)\cdot [u-(1+R)] 
         \\
         \Delta_0S(0)\cdot [d-(1+R)] 
       \end{cases}
       \end{align*}
   with respective probabilities $p$ and $1-p$. 
   
   If $(1+R)<d$, both quantities are positive and we could borrow money to buy assets
   to have a sure win. If $(1+R)>u$, both quantities are negative and we could 
   make money by selling assets and buying bonds. If $V(0)\neq0$, the argument is
   equally valid. 
   
   The arbitrage situation is not realistic if all the actors have complete 
   information. Thus, usually there is no free lunch!

\subsection{European options} 
    An \Conc{European call option} is a financial derivative: 
    It gives the holder the right (not the obligation) to buy a share for an  pre-specified amount 
    (\Conc{exercise price $K$})
    on a specific later date (\Conc{expiry date $T$}).     
    Similarly, an \Conc{European put option} is the right to sell the share. 
    
    If $S(T)$ is the value of the share at time $T$, the payoff of a call is
    $\big(S(T)-K\big)^+$. If $S(T)<K$, the holder does not exercise the option,
    since it can buy the share in the market for a cheaper price, so the payoff 
    is never negative.
    
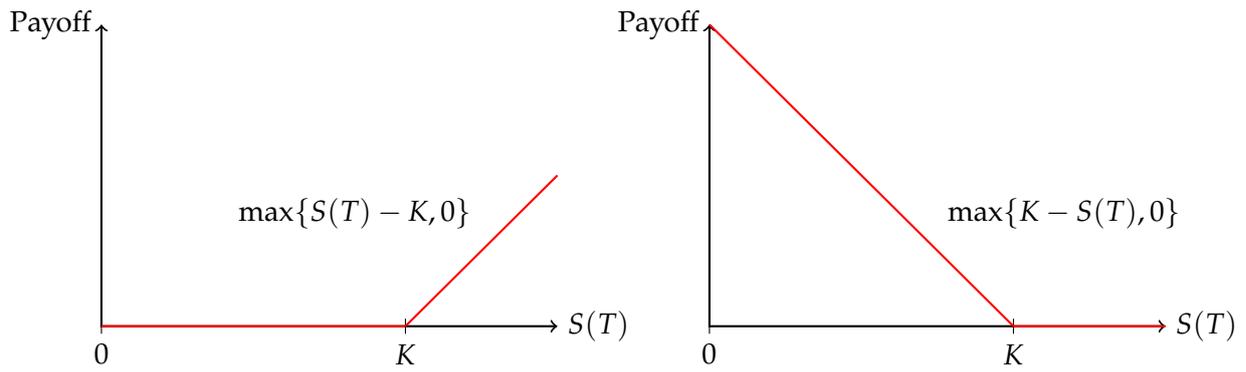
\begin{figure}
  \begin{center}
    \begin{tikzpicture}[scale=1] 
    \node [left] at (0,4) {Payoff};
    \node [right] at (6,0) {$S(T)$};
    \draw [<->, thick] (0,4) -- (0,0) -- (6,0);
    \draw [-] ( 0, -0.1) node [below] {$0$} -- ( 0, 0.1); 
    \draw [-] ( 4, -0.1) node [below] {$K$} -- ( 4, 0.1); 
    \draw [-, red, thick] ( 0, 0)  -- ( 4, 0);           
    \draw [-, red, thick] ( 4, 0)  -- ( 6, 2);           
    \node [left] at (5,1.5) {$\max\{S(T)-K,0\}$};

    \node [left] at (8,4) {Payoff};
    \node [right] at (14,0) {$S(T)$};
    \draw [<->, thick] (8,4) -- (8,0) -- (14,0);
    \draw [-] ( 8, -0.1) node [below] {$0$} -- ( 8, 0.1); 
    \draw [-] ( 12, -0.1) node [below] {$K$} -- ( 12, 0.1); 
    \draw [-, red, thick] ( 8, 4)  -- ( 12, 0);           
    \draw [-, red, thick] ( 12, 0)  -- ( 14, 0);           
    \node [right] at (11,1.5) {$\max\{K-S(T),0\}$};
    \end{tikzpicture}
  \end{center}
  \caption{Graphs of an European call and an European put}\label{callput}
\end{figure}

    Correspondingly, the payoff of a put is $\big(K-S(T)\big)^+$, see Figure \ref{callput}.

\subsection{Fair price of an European call option. Example} 
Assume the following data:
  \begin{itemize}
  \item 
    Current price of the share: $S(0)=100$
  \item 
    Interest of the risk-free bond: $10\%$ 
  \item 
    Possible prices for the share at time 1: $120$ or $90$
  \item 
    Exercise price: $K=100$ 
  \end{itemize}
    
    We have\quad $u=1.2$,\quad $d=0.9$,\quad $1+R=1.1$. 
    The payoff will be $C_u:=20$ or $C_d:=0$.
    
    To find the fair price, let us construct a portfolio with a value $V(1)$ equal to the
    payoff of the option. The fair price will be $V(0)$.    
    \begin{align*}
    V(1) &=
    \begin{cases}
    & B_0\cdot 1.1 + \Delta_0S(0)\cdot 1.2 = 20
    \\
    & B_0\cdot 1.1 + \Delta_0S(0)\cdot 0.9 = 0
    \end{cases}
    \\&\Rightarrow \Delta_0S(0) = 66.67 \ , \quad  B_0 = -54.55
    \end{align*}
    
    The fair price is thus 12.12

\subsection{Fair price of an European call option. In general} 
    In general, we have 
    \begin{align*}
    &V(1) =
    \begin{cases}
    &B_0(1+R) + \Delta_0S(0)u = C_u
    \\
    &B_0(1+R) + \Delta_0S(0)d = C_d
    \end{cases}
    \\
    &\Rightarrow
    B_0 = \frac{uC_d-dC_u}{(1+R)(u-d)} 
    \ ,\quad
    \Delta_0S(0)=\frac{C_u-C_d}{u-d}
    \ ,\quad
    \\
    &\Rightarrow
    B_0+\Delta_0S(0)=(1+R)^{-1}\big(C_uq+C_d(1-q)\big)
    \ ,
    \quad 
    \text{where $q=\frac{1+R-d}{u-d}$ .}
    \end{align*}
    
    It follows that the fair price of the option is the expected 
    (and discounted!) payoff of the option under the probability $Q=(q,1-q)$
    for the states of the economy:
  \begin{equation}\label{ExpQ}  
      \E_Q\big[(1+R)^{-1}\big(S(1)-K\big)^+\big]
  \end{equation}

  Some remarks on the probability $Q$: 
  \begin{itemize}
    \item
    Under $Q$, the share and the bond have the same expected return:
    \begin{align*}
    \E_Q[S(1)] &= S(0)uq+S(0)d(1-q)
    \\
    &=
    S(0)\Big(u\frac{1+R-d}{u-d}+ d\frac{u-1-R}{u-d}\Big)
    \\
    &= S(0)(1+R)
    \ .
  \end{align*}
    
  \item    
    The probability $Q$ does not depend on the underlying probability $P=(p,1-p)$ nor 
    on the payoff of the option.
    
  \item 
    $Q$ is called the \Conc{risk-neutral probability} (or \Conc{martingale probability}).
      \end{itemize}

\subsection{Fair price of an European call option. Example (cont.)}
  With the same data as before, we compute now the fair price directly using formula 
  (\ref{ExpQ}), where in this case $q=:2/3$. 
    \begin{align*}
    E_Q\big[(1+R)^{-1}\big(S(1)-K\big)^+\big] 
    &=(1+R)^{-1}\big(C_u\cdot q+C_d\cdot (1-q)\big)
    \\
    &=
    \frac{1}{1.1}\Big[20\cdot\frac{2}{3}+0\cdot\frac{1}{3}\Big]= 12.12 \ .
    \end{align*}
    
  Assume now that the exercise price is fixed to $K=95$ instead of $K=100$, 
  while all other data remain the same. Logically, the option should be more
  expensive in this case. Applying again formula (\ref{ExpQ}), 
    \begin{align*}
   E_Q\big[(1+R)^{-1}\big(S(1)-K\big)^+\big] 
    &=(1+R)^{-1}\big(C_u\cdot q+C_d\cdot (1-q)\big)
    \\
    &=
    \frac{1}{1.1}\Big[25\cdot\frac{2}{3}+0\cdot\frac{1}{3}\Big]= 15.15\ .
    \end{align*}

\subsection{European call option. Multiperiod} 
  The previous sections dealt with a single time period. Assume now that the 
  expiry time of the option is $T$ and that we can change the composition
  of the portfolio at any of the intermediate integer times. 
  
  A \Conc{trading 
  strategy} is then $\{(B_t,\Delta_t),\ 0\le t \le T-1\}$.
  It is called a \Conc{self-financing strategy} if we do not put new money
    or take money out of the portfolio.

      At time $t$, we can change the portfolio composition, but the value 
      remains the same:
    \begin{align*}
    B_{t}+\Delta_{t}S(t) = B_{t+1}+\Delta_{t+1}S(t) 
    \ .
    \end{align*}
      The new value at time $t+1$ will be:
    \begin{align*}
    B_{t+1}(1+R)+\Delta_{t+1}S(t+1) 
    \ .
    \end{align*}
      Therefore, the value increments for a self-financing strategy is 
    \begin{align*}
    V(t+1)-V(t) =B_{t+1}R+\Delta_{t+1}\big(S(t+1)-S(t)\big)
    \ .
    \end{align*}
    
  We can compute the fair price $F(0)$ at time 0 of an option with exercise value 
  $K$ at expiry date $T$ recursively: 
    \begin{align*}
    F (T-1) &=\E_Q\Big[\condprob{(1+R)^{-1}(S(T)-K)^+}{S(T-1)}\Big]    
    \\
    F(T-2) &=\E_Q\Big[\condprob{(1+R)^{-1}F(T-1)}{S(T-2)}\Big] 
    \\
    &=\E_Q\Big[\condprob{(1+R)^{-2}(S(T)-K)^+}{S(T-2)}\Big]
    \\
    &\phantom{=}\vdots
    \\
    F(0)&=\E_Q\big[(1+R)^{-T}(S(T)-K)^+\big] 
    \ .
    \end{align*}
  This computation uses essential properties of the conditional expectation that we are not
  going to detail here. But the conclusion must be quite intuitive.

\subsection{Martingales} 
    Under probability $Q$, the stochastic process $(1+R)^{-t}S(t)$ enjoys the martingale property.
    A stochastic process $\{X_t, t\ge 0\}$ is a \Conc{martingale} if
    \begin{align}\label{eq:martingale}
    \E\Big[\condprob{X_t}{X_s}\Big] = X_s  \quad\text{whenever $s<t$,}      
    \end{align}
    meaning that the knowledge of the state of the system at time $s$ makes this the expected 
    value at any later time. The discrete time process defined in Section \ref{sec:WrightFisher}
    is a discrete time martingale (see Equation (\ref{eq:dtmartingale})).   
    
    Martingales are good models for \emph{fair games}: The expected wealth of a player in the
    future is the current wealth, no matter what happened before, or how long has been playing.
    
    From (\ref{eq:martingale}) it can be deduced in particular that the expectation of the 
    process is constant in time. In our case of the European call option, this means
        \begin{align*}
    \E_Q\big[(1+R)^{-t}S(t)\big] = S(0)  
    \ ,
    \end{align*}
    implying that
    \begin{align*}
    \E_Q\big[S(T)\big] = S(0)\cdot (1+R)^T  
    \end{align*}
    which is precisely the return of the risk-free bond 
    (and this is why $Q$ is called a ``risk neutral'' probability measure).

\subsection{European call option. Continuous time} 
    In continuous time, 
        it can be shown that there is also a probability $Q$ under which $e^{-Rt}S(t)$ 
        is a martingale, and 
        the fair price at time 0 of a call option is given by
        \begin{align*}
        F(0)= \E_Q\big[e^{-RT}(S(T)-K)^+\big]
        \ ,
        \end{align*}
    although $Q$ is more difficult to describe here.
    
    The evolution of the value of the bond asset $I(t)$ is driven by the well-known
    differential equation 
    \begin{align*}
    dI(t)=R\cdot I(t)dt
    \end{align*}

    The evolution of the price of the share can be described as
    \begin{align}\label{eq:evolshare}
    dS(t)=S(t)\big(\mu dt + \sigma dW(t)\big)
    \end{align}
    where $W$ is a Wiener process, approximating (in the continuum limit) the 
    Markov chain given by the binomial model. The trend, if $p\neq 1/2$, goes
    to the \Conc{drift} $\mu$. The \Conc{volatility} $\sigma$ is the intensity
    of the noise.
    This is a simple example of a \Conc{stochastic differential equation}. 
    It is a pathwise description
    of a diffusion process with $b(t,x)\equiv\mu$ and $a(t,x)\equiv\sigma^2$ (see Section \ref{sec:diffusions}).

    Although the paths of $W$ are non-differentiable everywhere, Equation (\ref{eq:evolshare}) has the
    obvious meaning 
    \begin{align*}
    S(t)=S(0)+\mu\int_0^t S(r)\,dr + \sigma W(t) 
    \ .
    \end{align*}
    
    This equation 
    can be solved explicitly (this is not common, of course). The
    solution is the stochastic process given by 
    \begin{align*}
    S(t)=S(0)\exp\Big\{\mu t - \frac{1}{2}\sigma^2 t +\sigma W(t)\Big\}
    \ ,
    \end{align*}
    and we can compute its law from here.
    
    The evolution of the whole portfolio value will be
    \begin{align*}
    dV(t)&=B_t dI(t) + \Delta_t dS(t) 
    \ .
    \end{align*}

\subsection{Stochastic differential equations} 
  In general, a diffusion process $X$ with characteristic functions $a(t,x)$ and $b(t,x)$ 
  (called respectively \Conc{diffusion} and 
  \Conc{drift coefficients}) can be represented pathwise by means of the 
  stochastic differential equation 
\begin{equation*}
  dX(t)=b(t,X(t))\,dt+a(t,X(t))^{1/2}\,dW(t)
  \ ,
\end{equation*}  
  with a suitable definition of the last term, which in general, when the function $a$ depends
  effectively of its second argument, does not possess an obvious meaning.
 
 Diffusions can therefore be studied at the same time with the tools of partial differential
 equations that describe the evolution of the laws in time, and with the tools of stochastic
 processes and stochastic differential equations, that provide the evolution of the paths
 themselves.
 
 The word ``diffusion'' is taken from the physical phenomenon with that name: The movement of 
 particles in a fluid from regions of high concentration to regions of low concentration. 
 The heat ``diffuses'' in the same way, following the negative gradient of the temperature
 field $f(t,x)$. In one space dimension, it obeys the partial differential equation
\begin{equation*}
  \frac{\partial}{\partial t} f(t,x) = D\frac{\partial^2}{\partial x^2} f(t,x)
  \ , 
\end{equation*}  
  where $D$ is called the \emph{thermal diffusivity}. Comparing with Kolmogorov equations 
  (\ref{eq:Kolmof}-\ref{eq:Kolmob}),
  we see that, with suitable initial conditions, $f(t,x)$ is the density at time $t$ and point
  $x$ of a diffusion process following the stochastic differential equation
\begin{equation*}
  dX(t)=\sqrt{2D}\,dW(t)
\ , 
\end{equation*}  
  that means, essentially, the Wiener process. 
\section{Recommended books} 
  \begin{itemize}
  \item  
    Nelson, {\color{DarkGreen}Stochastic Modeling}, Dover 1995
    \newline
    (Arrivals, queues, Markov chains, simulation)
  \item  
    Gross-Harris, {\color{DarkGreen}Fundamentals of queueing theory}, Wiley 1998
    \newline
    (Queues)
  \item  
    Asmussen-Glynn, {\color{DarkGreen}Stochastic simulation}, Springer 2007
    \newline
    (Simulation)
  \item  
    Maruyama, {\color{DarkGreen}Stochastic problems in population genetics}, Springer 1977
    \newline
    (Diffusions, application to genetics)
  \item  
    Lamberton-Lapeyre, {\color{DarkGreen}Introduction au calcul stochastique appliqu\'{e} \`{a} la finance}, Ellipses 1997
    \newline
    (Diffusions, stochastic differential equations, application to finance)
  \end{itemize}

\newpage
\section*{Dedication: \Conc{Io non crollo}}

\begin{center}
  \includegraphics[scale=0.6]{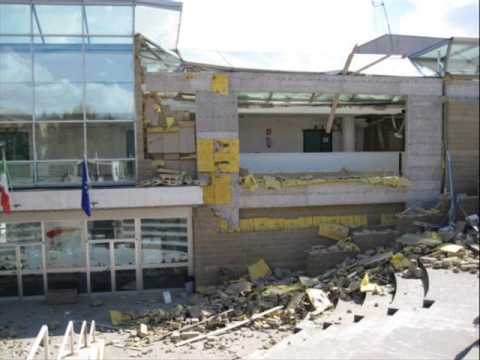}
\end{center}

This survey is based on a course given by the author in the Università degli Studi dell'Aquila, as a part
of the Intensive Programme \Conc{Mathematical Models in Life and Social Sciences},
in July 2008.

The year after, on April 6th 2009, the building where the programme took place was destroyed by
a strong earthquake that caused more than 300 deaths in the region.

This work is dedicated to the people that died, lost a beloved one, or lost their homes
that day. I adhere to the motto that helped the university people to carry on after the disaster: 

\centering\Conc{Io non crollo}.

\begin{center}
 \begin{overpic}
   [scale=0.5,tics=10]{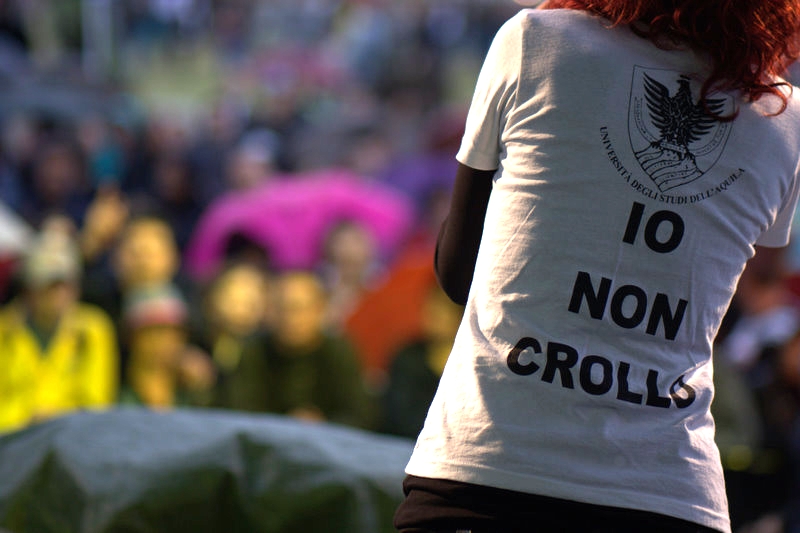}
   \put(55,-5) {\small Photo: Renato di Bartolomeo}
  \end{overpic}
\end{center}
\end{document}